\RequirePackage{fix-cm} 
\documentclass[a4paper, twoside, reqno, 12pt]{amsart}
\usepackage{fixltx2e}   

\usepackage{etex}

\usepackage[latin1]{inputenc}
\usepackage[T1]{fontenc}

\usepackage{esint}
\usepackage{dsfont}
\usepackage{xspace}
\usepackage{amsgen}
\usepackage{amsthm}
\usepackage{amssymb}
\usepackage{amsmath}
\usepackage{wasysym}
\usepackage{upgreek}
\usepackage{amsfonts}
\usepackage{stmaryrd}
\usepackage{mathtools}

\usepackage[mathcal, mathscr]{euscript}

\usepackage{mathrsfs}
\DeclareMathAlphabet{\mathscrbf}{OMS}{mdugm}{b}{n}

\usepackage{multicol}
\usepackage{multirow}

\usepackage{a4wide}

\usepackage[dvipsnames, table]{xcolor}
\definecolor{bckg}{RGB}{20.8, 20.8, 20.8}
\definecolor{oneblue}{rgb}{0.0, 0.0, 0.85}
\definecolor{Lightblue}{RGB}{214, 214, 214}
\definecolor{bluepigment}{rgb}{0.2, 0.2, 0.6}
\definecolor{charcoal}{rgb}{0.21, 0.27, 0.31}
\definecolor{denimblue}{rgb}{0.08, 0.38, 0.74}
\definecolor{Lightgray}{rgb}{0.89, 0.89, 0.89}
\definecolor{darkgrey}{rgb}{0.273, 0.281, 0.30}
\definecolor{darkelectricblue}{rgb}{0.33, 0.41, 0.47}

\usepackage[sort&compress, comma, square, numbers]{natbib}

\usepackage{psfrag}
\usepackage{graphicx}
\usepackage{subfigure}
\usepackage{morefloats}
\usepackage{indentfirst}

\usepackage{acronym}
\usepackage{microtype}
\usepackage[labelsep=period,%
            labelfont={bf,sf,color=bluepigment},%
            justification=raggedright]{caption}

\usepackage[perpage, symbol]{footmisc}

\usepackage[usenames, pdf]{pstricks}
\usepackage{epsfig}
\usepackage{pst-grad} 
\usepackage{pst-plot} 

\usepackage{url}
\usepackage[colorlinks,
           urlcolor=oneblue,
           linkcolor=denimblue,
           citecolor=NavyBlue,
           bookmarksopen=false,
           pdfpagemode=UseNone,
           pagebackref]{hyperref}

\usepackage[explicit]{titlesec}

\titleformat{\section}[block]
  {\color{NavyBlue}\Large\sffamily\bfseries}
  {}
  {0.0em}
  {\colorbox{bckg!5}{\strut\parbox{\dimexpr\linewidth-2\fboxsep\relax}{\thesection. #1}}}
  [\vspace*{0.33em}]

\titleformat{name=\section,numberless}[block]
  {\color{NavyBlue}\Large\sffamily\bfseries}
  {}
  {0.0em}
  {\colorbox{bckg!5}{\strut\parbox{\dimexpr\linewidth-2\fboxsep\relax}{#1}}}
  [\vspace*{0.33em}]

\titleformat{\subsection}
  {\color{NavyBlue}\large\sffamily\bfseries}
  {}
  {0.0em}
  {\colorbox{bckg!5}{\parbox{\dimexpr\linewidth-2\fboxsep\relax}{\thesubsection. #1}}}
  [\vspace*{0.33em}]

\titleformat{name=\subsection,numberless}
  {\color{NavyBlue}\Large\sffamily\bfseries}
  {}
  {0em}
  {\colorbox{bckg!5}{\parbox{\dimexpr\linewidth-2\fboxsep\relax}{#1}}}
  [\vspace*{0.33em}]

\titleformat{\subsubsection}
  {\color{bluepigment}\sffamily\normalsize\bfseries}
  {\thesubsubsection}
  {0.5em}
  {#1}
  [\vspace*{0.33em}]

\titleformat{\paragraph}[runin]
  {\color{bluepigment}\sffamily\small\bfseries}
  {}
  {0em}
  {#1}

\titlespacing{\section}{0.0em}{1.5em plus 2pt minus 2pt}%
{1.0em plus 2pt minus 2pt}[0em]
\titlespacing{\subsection}{0.5em}{1.5em plus 2pt minus 2pt}%
{1.0em}[0em]
\titlespacing{\subsubsection}{0.5em}{1.5em plus 2pt minus 2pt}%
{1.0em plus 2pt minus 2pt}[0em]

\usepackage{titletoc}

\contentsmargin{0.5em}
\newlength{\tocsep} 
\titlecontents{section}[\tocsep]
  {\addvspace{10pt}\bfseries\sffamily}
  {\contentslabel[\thecontentslabel]{\tocsep}}
  {}
  {\ \titlerule*[0.75pc]{.}\ \thecontentspage}
  []
\titlecontents{subsection}[\tocsep]
  {\addvspace{8pt}\sffamily}
  {\contentslabel[\thecontentslabel]{\tocsep}}
  {}
  {\ \titlerule*[0.5pc]{.}\ \thecontentspage}
  []
\titlecontents*{subsubsection}[\tocsep]
  {\addvspace{2pt}\footnotesize\sffamily}
  {}
  {}
  {\ \titlerule*[0.35pc]{.}\ \thecontentspage}
  [\\*]

\makeatletter
\def\@setauthors{%
  \begingroup
  \def\thanks{\protect\thanks@warning}%
  \trivlist
  \centering\footnotesize \@topsep30\p@\relax
  \advance\@topsep by -\baselineskip
  \item\relax
  \author@andify\authors
  \def\\{\protect\linebreak}%
  \textsc{\normalsize\textcolor{darkelectricblue}{\authors}}%
  \ifx\@empty\contribs
  \else
    ,\penalty-3 \space \@setcontribs
    \@closetoccontribs
  \fi
  \endtrivlist
  \endgroup
}
\def\@settitle{\begin{center}%
  \baselineskip14\p@\relax
    \bfseries
    \textsc{\Large\textcolor{charcoal}{\@title}}
  \end{center}%
}
\makeatother

\usepackage{enumitem}
\setlist[description]{%
  topsep=30pt,               
  itemsep=5pt,               
  font={\bfseries\sffamily\color{NavyBlue}}, 
}

\usepackage{fancyhdr}
\usepackage{lastpage}

\newcommand*\Title{\textcolor{bluepigment}{On the estimation of moisture permeability and advection coefficients}}
\newcommand*\Authors{\textcolor{bluepigment}{J.~Berger, T.~Busser, D.~Dutykh \& N.~Mendes}}
\newcommand*{\plogo}{\textcolor{gray}{{\texttt{arXiv.org} / \textsc{hal}}}} 

\numberwithin{equation}{section}

\newcommand*\unit[1]{\mathsf{#1}}
\newcommand*\unitfrac[2]{\frac{\mathsf{#1}}{\mathsf{#2}}}

\newcommand*\pd[2]{\frac{\partial #1}{\partial #2}}

\newcommand{\eqdef}{\mathop{\stackrel{\,\mathrm{def}}{:=}\,}}

\renewcommand{\O}{\mathcal{O}}

\newcommand*\e[1]{\cdot 10^{\,#1}}

\newcommand*\egal{\ = \ }
\newcommand*\plus{\ + \ }
\newcommand*\moins{\ - \ }

\newcommand{\dPhi}{\Delta \phi_{\,\infty}}

\newcommand*\vP{\mathbf{P}}

\newcommand*{\Ox}{\Omega_{\, x}}

\newcommand{\f}{\mathrm{f}}

\newcommand{\kl}{k_{\,l}}
\newcommand{\kv}{k_{\,v}}
\newcommand{\Pc}{P_{\,c}}
\newcommand{\Ps}{P_{\,s}}
\newcommand{\Pv}{P_{\,v}}
\newcommand{\Pvi}{P_{\,v}^{\,i}}

\newcommand{\Rv}{R_{\,v}}
\newcommand{\TL}{T^{\,\infty}}
\newcommand{\rholv}{\rho_{\,l+v}}

\newcommand{\phiL}{\phi^{\,\infty}}

\newcommand{\eg}{\emph{e.g.}\xspace}
\newcommand{\etc}{\emph{etc.}\xspace}

\begin{document}

\title[\Title]{On the estimation of moisture permeability and advection coefficients of a wood fibre material using the optimal experiment design approach}

\author[J.~Berger]{Julien Berger$^*$}
\address{\textbf{J.~Berger:} LOCIE, UMR 5271 CNRS, Universit\'e Savoie Mont Blanc, Campus Scientifique, F-73376 Le Bourget-du-Lac Cedex, France}
\email{Berger.Julien@hotmail.com}
\urladdr{https://www.researchgate.net/profile/Julien\_Berger3/}
\thanks{$^*$ Corresponding author}

\author[T.~Busser]{Thomas Busser}
\address{\textbf{T.~Busser:} LOCIE, UMR 5271 CNRS, Universit\'e Savoie Mont Blanc, Campus Scientifique, F-73376 Le Bourget-du-Lac Cedex, France}
\email{Thomas.Busser@univ-smb.fr}
\urladdr{https://www.researchgate.net/profile/Thomas\_Busser/}

\author[D.~Dutykh]{Denys Dutykh}
\address{\textbf{D.~Dutykh:} LAMA, UMR 5127 CNRS, Universit\'e Savoie Mont Blanc, Campus Scientifique, F-73376 Le Bourget-du-Lac Cedex, France}
\email{Denys.Dutykh@univ-savoie.fr}
\urladdr{http://www.denys-dutykh.com/}

\author[N.~Mendes]{Nathan Mendes}
\address{\textbf{N.~Mendes:} Thermal Systems Laboratory, Mechanical Engineering Graduate Program, Pontifical Catholic University of Paran\'a, Rua Imaculada Concei\c{c}\~{a}o, 1155, CEP: 80215-901, Curitiba -- Paran\'a, Brazil}
\email{Nathan.Mendes@pucpr.edu.br}
\urladdr{https://www.researchgate.net/profile/Nathan\_Mendes/}

\keywords{inverse problem; parameter estimation; optimal experiment design (OED); convective moisture transport; sensitivity functions; model identification}

\begin{titlepage}
\thispagestyle{empty} 
\noindent
{\Large Julien \textsc{Berger}}\\
{\it\textcolor{gray}{Polytech Annecy--Chamb\'ery, LOCIE, France}}
\\[0.02\textheight]
{\Large Thomas \textsc{Busser}}\\
{\it\textcolor{gray}{Polytech Annecy--Chamb\'ery, LOCIE, France}}
\\[0.02\textheight]
{\Large Denys \textsc{Dutykh}}\\
{\it\textcolor{gray}{CNRS--LAMA, Universit\'e Savoie Mont Blanc, France}}
\\[0.02\textheight]
{\Large Nathan \textsc{Mendes}}\\
{\it\textcolor{gray}{Pontifical Catholic University of Paran\'a, Brazil}}
\\[0.08\textheight]

\vspace*{0.7cm}

\colorbox{Lightblue}{
  \parbox[t]{1.0\textwidth}{
    \centering\huge\sc
    \vspace*{0.7cm}
    
    \textcolor{bluepigment}{On the estimation of moisture permeability and advection coefficients of a wood fibre material using the optimal experiment design approach}
    
    \vspace*{0.7cm}
  }
}

\vfill 

\raggedleft     
{\large \plogo} 
\end{titlepage}


\newpage
\thispagestyle{empty} 
\par\vspace*{\fill}   
\begin{flushright} 
{\textcolor{denimblue}{\textsc{Last modified:}} \today}
\end{flushright}


\newpage
\maketitle
\thispagestyle{empty}


\begin{abstract}

This paper presents a practical application of the concept of Optimal Experiment Design (OED) for the determination of properties of porous materials with \emph{in-situ} measurements and an identification method. First, an experimental set-up was presented and used for the measurement of relative humidity within a wood fibre material submitted to single and multiple steps of relative humidity variation. Then, the application of OED enabled to plan the experimental conditions in terms of sensor positioning and  boundary conditions out of 20 possible designs. The OED search was performed using the Fisher information matrix and \emph{a priori} knowledge of the parameters. It ensures to provide the best accuracy of the identification method and thus the estimated parameter. Optimal design results have been found for single steps from the relative humidity $\phi \egal 10$ to $75 \ \unit{\%}\,$, with one sensor located at the position $X$ between $4$ and $6 \ \unit{cm}\,$, for the estimation of moisture permeability coefficients, while from $\phi \egal 75 \ \unit{\%}$ to $\phi \egal 33 \ \unit{\%}\,$, with one sensor located at $ X^{\,\circ} \egal 3 \ \unit{cm} \,$, for the estimation of the advection coefficient. The OED has also been applied for the identification of couples of parameters. A sample submitted to multiple relative humidity steps ($\phi \egal 10-75-33-75 \ \unit{\%}$) with a sensor placed at $ X^{\,\circ} \egal 5 \ \unit{cm}$ was found as the best option for determining both properties with the same experiment. These OED parameters have then been used for the determination of moisture permeability and advection coefficients. The estimated moisture permeability coefficients are twice higher than the \emph{a priori} values obtained using standard methods. The advection parameter corresponds to the mass average velocity of the order of $\mathsf{v} \egal 0.01 \ \unit{mm/s}$ within the material and may play an important role on the simulation of moisture front.

\bigskip
\noindent \textbf{\keywordsname:} inverse problem; parameter estimation; optimal experiment design (OED); convective moisture transport; sensitivity functions; model identification \\

\smallskip
\noindent \textbf{MSC:} \subjclass[2010]{ 35R30 (primary), 35K05, 80A20, 65M32 (secondary)}
\smallskip \\
\noindent \textbf{PACS:} \subjclass[2010]{ 44.05.+e (primary), 44.10.+i, 02.60.Cb, 02.70.Bf (secondary)}

\end{abstract}


\newpage
\tableofcontents
\thispagestyle{empty}


\newpage
\section{Introduction}

Heating or cooling strategies and design of building envelopes are commonly based on numerical simulations performed by hygrothermal tools such as Delphin \cite{BauklimatikDresden2011}, MATCH \cite{Rode2003}, MOIST \cite{Burch1993}, WUFI \cite{IBP2005} or Umidus \cite{Mendes1997, Mendes1999} or by whole-building simulation programs, used in the frame of the International Energy Agency Annex $41$ \cite{Woloszyn2008}. Details about the mathematical models and their numerical schemes are described in \cite{Mendes2017}. Those programs are capable to simulate whole buildings considering the combined heat and moisture transfer through porous elements.

Nevertheless, porous materials have hygrothermal properties strongly dependent on moisture content --- mainly highly hygroscopic materials such as wood, are commonly estimated using experimental characterisation methods. The moisture capacity is traditionally determined using the gravimetric methods (\textsf{ISO 12571}). Samples are weighed for different relative humidity conditions, ensured by salt solutions. For the vapour permeability, the most common measurement procedure is based on the standard (\textsf{ISO 12572}). It measures the mass variation of a material sample under a controlled difference of relative humidity between both sides. The vapour permeability is estimated for different relative humidity levels (dry and wet cup). For the liquid permeability, the sample is exposed to a liquid pressure gradient. The permeability is estimated by measuring the liquid flux at the equilibrium.

Even if measurements of these properties are well established and largely reported in the literature, some discrepancies can appear when comparing experimental data from \emph{in-situ} measurements to numerical model results. A material, with an initial moisture content $w_{\,0}\,$, is submitted to an adsorption phase at relative humidity $\phi_{\,1}$ and then to a desorption phase at relative humidity $\phi_{\,2}\,$. Results of the simulation underestimate the adsorption process or overestimate the desorption process. Numerous studies state similar observations. Interested readers may consult \cite{James2010, Berger2017a} for a preliminary introduction to this investigation.

To answer this issue, models are calibrated using \emph{in-situ} measurements for estimating material properties to reduce the discrepancies between model predictions and real performance. For instance, in \cite{Colinart2016b}, liquid and vapour transfer coefficients are determined for hemp material submitted to a drying process. In \cite{Kanevce2005}, moisture- and temperature-dependent diffusivity and thermo-physical properties are estimated using only temperature measurements under a drying process. In \cite{Rouchier2017}, moisture properties of a wood fibre material were estimated and so the parameter estimation problem revealed a resistance factor for vapour permeability lower than $1$, which is physically unacceptable as the vapour diffusion of through the porous material cannot be faster than in the air. Physical models commonly take into account only diffusion transfer and authors conclude that other physical phenomena, such as moisture advection, might not be ignored. Therefore, this research has been conducted to estimate the moisture permeability and advection coefficients of a similar wood fibre material. First, an experimental set-up has been used, enabling to submit the material samples to several steps of relative humidity, under isothermal conditions. Temperature and relative humidity sensors have been placed within the samples.

The estimation of the unknown parameters, \eg, wall thermo-physical properties, based on observed data and identification methods, strongly depends on the experimental protocol and particularly on the imposed boundary conditions and on the location of the sensors. In \cite{Berger2017}, the concept of searching the Optimal Experiment Design (OED) was used to determine the best experimental conditions in terms of quantity and location of sensors, and flux imposed to the material. These conditions ensure to provide the best accuracy of the identification method and thus the estimated parameter. It was studied for the identification of hygrothermal properties of a porous material. Even though the approach was verified by solving $100$ inverse problems for different experiment designs, it remained theoretical as no real experimental data nor existing facility were considered.

This research intends to go one step further. The issue is to use the methodology to determine the OED according to the existing facility. It first aims at defining the optimal boundary condition and location of the sensor to ensure an accurate solution of the parameter estimation problem. Then, experimental measurements are provided respecting the particular design and used to estimate the unknown parameters. Thus, this article is organised as follows. Next Section~\ref{sec:methodology} presents the physical problem and the methodology of OED searching. Then, Section~\ref{sec:exp_facility} describes the experimental facility used for providing measurements in a wood fibre material. In Section~\ref{sec:OED}, the OED is searched according to the different possibilities for the estimation of one or several parameters. Section~\ref{sec:estim_parameters} provides an estimation of the moisture permeability and advection coefficients, given the experimental data obtained according to the OED.


\section{Methodology}
\label{sec:methodology}

\subsection{Physical problem and mathematical formulation}
\label{sec:moisture_equations}

\begin{figure}
\begin{center}
  \includegraphics[width=.8\textwidth]{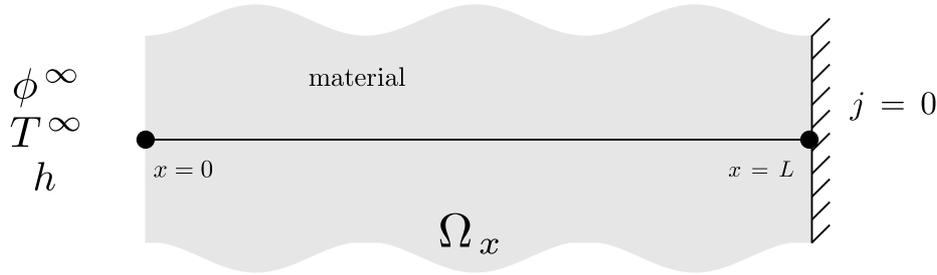}
  \caption{\small\em Illustration of the physical problem.}
  \label{fig:material}
\end{center}
\end{figure}

The physical problem involves one-dimension moisture convection through a porous material defined by the spatial domain $\Ox \egal [\, 0, \, L \,]$ as shown in Figure~\ref{fig:material}. The moisture transfer occurs due to capillary migration, vapour diffusion and advection of the vapour phase. The physical problem can be formulated as the convective moisture equation \cite{Luikov1966, Tariku2010, Belleudy2016}:
\begin{align}\label{eq:moisture_equation_1D}
  & \pd{\rholv}{t} \egal \pd{}{x} \left( \, \kl \, \pd{\Pc}{x} \plus \kv \, \pd{\Pv}{x} \, \right) \moins \pd{}{x}\left( \, \frac{\Pv}{\Rv \ T} \, \mathsf{v} \, \right) \,,
\end{align}
where $\rholv$ is the volumetric moisture content of the material, $\kv$ and $\kl\,$, the vapour and liquid permeabilities, $\Pv\,$, the vapour pressure, $\Pc\,$, the capillary pressure, $T\,$, the temperature, $\mathsf{v}\,$, the air velocity and, $\Rv\,$, the water vapour gas constant. Eq.~\eqref{eq:moisture_equation_1D} can be written using the vapour pressure $\Pv$ as the driving potential. For this, we consider the physical relation, known as the \textsc{Kelvin} equation, between $\Pv$ and $\Pc\,$:
\begin{align*}
  \Pc & \egal \Rv \, T \, \ln\left(\frac{\Pv}{\Ps(T)}\right) \,,\\
  \pd{\Pc}{\Pv} & \egal \frac{R_{\,v} \, T}{\Pv} \,.
\end{align*}
Thus we have:
\begin{align*}
  \pd{\Pc}{x} \egal \pd{\Pc}{\Pv} \, \pd{\Pv}{x} \plus \pd{\Pc}{T} \, \pd{T}{x} \,.
\end{align*}

The temperature remains the same at the boundaries. Even if heat transfer occurs in the material due to phase change, the temperature variations in the material are assumed negligible. Thus, the second right-hand side term vanishes and we obtain:
\begin{align*}
& \pd{\Pc}{x} \egal  \frac{\Rv \, T}{\Pv} \, \pd{\Pv}{x} \,.
\end{align*}
In addition, we have:
\begin{align*}
& \pd{\rholv}{t} \egal \pd{\rholv}{\phi} \, \pd{\phi}{\Pv} \, \pd{\Pv}{t} \plus \pd{\rholv}{T} \, \pd{T}{t} \,.
\end{align*}
Under isothermal conditions, the second right-hand term of the equation above also vanishes. Considering the relation $\rholv \egal \f \, (\phi) \egal \f \, (\Pv,\,T)\,$, obtained from material properties and from the relation between the vapour pressure $\Pv$ and the relative humidity $\phi\,$, we get:
\begin{align*}
 & \pd{\rholv}{t} \egal \frac{1}{\Ps} \; \pd{\rholv}{\phi} \; \pd{\Pv}{t} 
 \egal \frac{1}{\Ps} \; \pd{\f}{\phi} \; \pd{\Pv}{t} \,. 
\end{align*}

The advection flow is given by:
\begin{align*}
  j_{\,a} & \egal \egal \rho_{\,k} \, \mathsf{v} \,,
\end{align*}
where $\mathsf{v}$ is the molar average velocity and $\rho_{\,k}$ is the volumetric concentration of moving moisture mass, calculated as: 
\begin{align*}
  \rho_{\,k} \egal \frac{\rho_{\,m}}{b} \,,
\end{align*}
where $\rho_{\,m}$ is the volumetric concentration of moisture and $b$ is the ratio of the volume of moisture $V_{\,m}$ to the total volume of capillaries $V$. Thus, we have: 
\begin{align*}
  b \ \eqdef \ \frac{V_{\,m}}{V}\,.
\end{align*}
It is also assumed that there is no variation of the volume of capillaries. In other words, the shrinkage or expansion effects, due to the variation of the moisture content, are neglected \footnote{This hypothesis could be revised by writing the volume of capillaries as $V \egal \bigl(\, 1 \plus \beta \,\bigr) \, V_{\,0}\,$, where $V_{\,0}$ volume of capillaries at dry state and $\beta\,$, the shrinkage coefficient, depending on the moisture content $\rholv \,$.} Using the perfect gas law, we obtain: 
\begin{align*}
  j_{\,a} \egal \frac{\Pv}{R_{\,v} \, T} \; \mathsf{v} \,.
\end{align*}
With the assumption of isothermal conditions and constant mass average velocity $\mathsf{v}\,$, it can be written: 
\begin{align*}
  & \pd{}{x}\left( \, \frac{\Pv}{\Rv \ T} \, \mathsf{v} \, \right) \simeq \frac{\mathsf{v}}{\Rv \ T} \, \pd{\Pv}{x} \,.
\end{align*}
Eq.~\eqref{eq:moisture_equation_1D} can be therefore rewritten as:
\begin{align}\label{eq:moisture_equation_1D_v2}
  & \f^{\,\prime}(\Pv) \; \frac{1}{\Ps} \; \pd{\Pv}{t} \egal \pd{}{x} \Biggl[ \, \Bigl( \, \kl \, \frac{\Rv \, T}{\Pv} \plus \kv \, \Bigr) \, \pd{\Pv}{x} \, \Biggr] \moins \frac{\mathsf{v}}{\Rv \ T} \, \pd{\Pv}{x} \,.
\end{align}

The material properties $\f(\Phi)\,$, $\kl$ and $\kv$ depend on the vapour pressure $\Pv\,$. We denote $d \egal \kl \, \dfrac{\Rv \, T}{\Pv} \plus \kv $ the global moisture transport coefficient, $c \egal \dfrac{\partial \f}{\partial \phi} \; \dfrac{1}{\Ps}$ the moisture storage coefficient and, as $a \egal \dfrac{\mathsf{v}}{\Rv \ T}\,$, the advection coefficient. Thus, Eq.~\eqref{eq:moisture_equation_1D_v2} becomes:
\begin{align}\label{eq:HAM_equation2}
& c \; \pd{\Pv}{t} \egal \pd{}{x} \Biggl[ \, d \, \pd{\Pv}{x} \, \Biggr] \moins a \, \pd{\Pv}{x} \,.
\end{align}

At $x \egal 0\,$, the surface is in contact with the ambient air of temperature $\TL$ and relative humidity $\phiL$. Thus, the boundary condition is expressed as:
\begin{align}\label{eq:HAM_BC_L}
  d \, \pd{\Pv}{x} \moins a \, \Pv & 
  \egal h \, \biggl(\, \Pv - \Ps \bigl(\, \TL \,\bigr) \ \phiL \,\biggr) \,.
\end{align}
At $x \egal L\,$, the flow is imposed to be null: 
\begin{align}\label{eq:HAM_BC_R}
  d \, \pd{\Pv}{x} \moins a \, \Pv & \egal 0 \,.
\end{align}
At $t \egal 0\,$, the vapour pressure is supposed to be uniform in the material
\begin{align}\label{eq:HAM_ic}
 \Pv &\egal \Pvi \,.
\end{align}

The following assumptions are adopted on the properties of the material. The moisture capacity $c$ is assumed as a second-degree polynomial of the relative humidity. The moisture permeability $d$ is taken as a first-degree polynomial of the relative humidity:
\begin{align}\label{eq:mat_hypothesis}
  & c \egal c_{\,0} \plus c_{\,1} \, \phi \plus c_{\,2} \, \phi^{\,2} \,, \\
  & d \egal d_{\,0} \plus d_{\,1} \, \phi \,.
\end{align}
As the sorption isotherm function is a monotonic increasing function, a constraint is imposed on $c_{\,0}\,$, $c_{\,1}$ and $c_{\,2}\,$. Other functions can be used to describe the variation of the moisture permeability, for instance, exponential ones. However, the permeability of the material under investigation has been determined using the cup method in \cite{Rafidiarison2015, Vololonirina2014} and expressed as a first-order polynomial, which will be used as \emph{a priori} parameters in the next Sections.

The direct problem, defined by Eqs.~\eqref{eq:HAM_equation2}, \eqref{eq:HAM_BC_L}, \eqref{eq:HAM_BC_R} and \eqref{eq:HAM_ic}, is solved using a finite-difference standard discretisation. An embedded adaptive in time \textsc{Runge}--\textsc{Kutta} scheme, combined with a \textsc{Scharfetter}--\textsc{Gummel}  spatial discretisation approach, is used \cite{Berger2017a}. It is adaptive and embedded to estimate local error in time with little extra cost. The algorithm was implemented in the \texttt{Matlab} environment. It is also important to note that a unitless formulation of the problem was used, due to a number of good reasons. First of all, it enables to determine important scaling parameters (\textsc{Biot}, \textsc{Luikov} numbers for instance). Henceforth, solving one dimensionless problem is equivalent to solve a whole class of dimensional problems sharing the same scaling parameters. Then, dimensionless equations allow to estimate the relative magnitude of various terms, and thus, eventually to simplify the problem using asymptotic methods \citep{Nayfeh2000}. Finally, the floating point arithmetics is designed such as the rounding errors are minimal if you manipulate the numbers of the same magnitude \citep{Kahan1979}. Moreover, the floating point numbers have the highest density in the interval $(\, 0,\,1 \,)$ and their density decays exponentially when we move further away from zero. So, it is always better to manipulate numerically the quantities of the order of $\O(1)$ to avoid severe round-off errors and to likely improve the conditioning of the problem in hands.

\subsection{Optimal Experiment Design}
\label{subsec:OED}

Efficient computational algorithms for recovering parameters $\vP$ given an observation $u_{\, \mathrm{exp}}$ of the field $u\, (x,\, t)$ have already been proposed. Readers may refer to \cite{Ozisik2000} for a primary overview of different methods. They are based on the minimisation of the cost function $\mathrm{J}\, [ \, \vP \, ]\,$. For this, it is required to equate to zero the derivatives of $\mathrm{J}\, [ \, \vP \, ] \,$, with respect to each of the unknown parameters $p_{\,m}$ to find critical points. Associated to this necessary condition for the minimisation of $\mathrm{J}\, [ \, \vP \, ]\,$, the scaled dimensionless local sensitivity function \cite{Finsterle2015} is introduced:
\begin{align}\label{sec1_eq:sensitivity_matrix}
& \Theta_{\,m} \,(x,t)= \frac{\sigma_{\,p}}{\sigma_{\,u}} \,\pd{u}{p_{\,m}} \,, && \forall m \ \in \   \bigl\{ 1, \ldots, M \, \bigr\}
\end{align}
where $\sigma_{\,u}$ is the variance of the error measuring $u_{\, \mathrm{exp}}$. The parameter scaling factor $\sigma_{\,p}$ equals $1$ as we consider that prior information on parameter $p_{\,m}$ has low accuracy. It is important to note that all the algorithms has been developed considering the dimensionless problem in order to compare only the order of variation of parameters and observation, avoiding the effects of units and scales.

The function $\Theta_{\,m}$ measures the sensitivity of the estimated field $u$ with respect to change in the parameter $p_{\,m}$ \cite{Nenarokomov2005, Ozisik2000, Artyukhin1985}. A small magnitude of $\Theta_{\,m}$ indicates that large changes in $p_{\,m}$ induce small changes in $u\,$. The estimation of parameter $p_{\,m}$ is therefore difficult in this case. When the sensitivity coefficient $\Theta_m$ is small, the inverse problem is ill-conditioned. If the sensitivity coefficients are linearly dependent, the inverse problem is also ill-conditioned. Therefore, to get an optimal evaluation of parameters $\vP\,$, it is desirable to have linearly-independent sensitivity functions $\Theta_{\,m}$ with large magnitude  for all parameters $p_{\,m}\,$. These conditions ensure the best conditioning of the computational algorithm to solve the inverse problem and thus the better accuracy of the estimated parameter.

It is possible to define the experimental design in order to reach these conditions. The issue is to find the optimal sensor location $X^{\,\circ}$ and the optimal amplitude $\phi^{\,\infty, \,\circ}$ of the relative humidity of the ambient air at the material bounding surface, $x \egal 0\,$. To search this optimal experiment design, we introduce the following measurement plan:
\begin{align}\label{sec1_eq:measurement_plan}
  & \pi \egal \bigl\{ \, X \,, \phi^{\,\infty} \, \bigr\} \,.
\end{align}
Here, the measurement plan considers the whole values of the space domain $\Omega_{\,x}$ for the sensor position $X\,$. When solving numerically the problem, a step $\Delta X \egal 0.01\,$, between each sensor position, is taken into account. For the relative humidity $\phi^{\,\infty} \,$, the measurement plan is fixed according to the experimental facility. As detailed in the next sections, $20$ different boundary conditions are considered.

In the analysis of optimal experiments for estimating the unknown parameter(s) $\vP$, a quality index describing the recovering accuracy is the $D-$optimum criterion \cite{Nenarokomov2005, Karalashvili2015, Artyukhin1985, Beck1977, Ucinski2004, VandeWouwer2000, Sun2005, Fadale1995, Emery1998, Anderson2005, Alifanov1995, Terejanu2012}:
\begin{align}
\label{sec1_eq:D_optimum}
\Psi \egal \det \bigl[ \, F(\,\pi\,) \, \bigr] \,,
\end{align}
where $F(\,\pi\,)$ is the normalized \textsc{Fisher} information matrix \cite{Karalashvili2015, Ucinski2004} defined as:
\begin{subequations}\label{sec1_eq:fisher_matrix}
\begin{align}
   F(\,\pi\,) & \egal \bigl[\, \Phi_{\,i \, j} \,\bigr] \,, && \forall (i,j) \ \in \ \bigl\{ 1, \ldots, M \, \bigr\}^{\,2}\,, \\[3pt]
   \Phi_{\,i \, j} & \egal \sum_{\,n \,= \,1}^{\,N} \ \int_{\,0}^{\, \tau} \Theta_{\,i} \ (x_{\,n} \,, t) \ \Theta_{\,j} \ (x_{\,n} \,, t) \ \mathrm{dt} \,.
\end{align}
\end{subequations}
where $N$ is the number of sensors and $\tau$ the simulation horizon time.

The matrix $F(\,\pi\,)$ characterises the total sensitivity of the system as a function of measurement plan $\pi$ (Eq.~ \eqref{sec1_eq:measurement_plan}). The OED search aims at finding a measurement plan $\pi^{\,\circ}$ for which the objective function -Eq.~\eqref{sec1_eq:D_optimum}- reaches the maximum value:
\begin{align}\label{sec1_eq:optimal_experimental_design}
  & \pi^{\,\circ} \egal \bigl\{ \, X^{\,\circ} \,, \phi^{\,\infty, \,\circ} \, \bigr\} \egal \arg \max_{\,\pi} \Psi \,.
\end{align}

To solve the Eq.~\eqref{sec1_eq:optimal_experimental_design}, a domain of variation $\Omega_{\,\pi}$ is considered for the sensor position $X$ and the amplitude $\phi^{\,\infty}$ of the boundary conditions. Then, the following steps are carried out for each value of the measurement plan $\pi \egal \bigl\{ \, X \,, \phi^{\,\infty} \, \bigr\} $ in the domain $\Omega_{\,\pi}\,$. First, the direct problem, defined by Eqs.~\eqref{eq:HAM_equation2}-\eqref{eq:HAM_ic}, is solved. Then, given the solution $\Pv$ for a fixed value of the measurement plan, the next step consists of computing the sensitivity coefficients $\Theta_{\,m} \egal \dfrac{\partial u}{\partial p_{\,m}}\,$, using also an embedded adaptive time \textsc{Runge}--\textsc{Kutta} scheme combined with central spatial discretisation. Then, with the sensitivity coefficients, the \textsc{Fisher} matrix \eqref{sec1_eq:fisher_matrix}(a,b) and the $D-$optimum criterion \eqref{sec1_eq:D_optimum} are calculated. All in all, four partial differential equations have to be solved for each values of the measurement plan: one for the vapor pressure and three for the sensitivity functions associated to the three parameters under investigation. Using the \texttt{Matlab} platform and the \texttt{ODE45} function, the computation of the results for one value of the measurement plan takes $ 150 \ \unit{s}$ with Intel i7 CPU and 16GB of RAM. Thus, the exploration of the whole measurement plan requires less than $1 \ \unit{h}$. The solution of the direct and sensitivity problems are obtained for a given \emph{a priori} parameter $\vP$ and, in this case, the validity of the OED depends on this knowledge. If there is no prior information, the methodology of the OED can be done using an outer loop on the parameter $\vP$ sampled using, for instance, Latin hypercube or \textsc{Halton} or \textsc{Sobol} quasi-random samplings. Interested readers may refer to \cite{Berger2017} for further details on the computation of sensitivity coefficients.

Equation~\eqref{sec1_eq:measurement_plan} has been written for a relative humidity step at the boundary of the material. However, the methodology can be applied for any other boundary condition. For instance, in \cite{Berger2017}, time varying boundary conditions in relative humidity and temperature are considered for the analysis of the OED. In \cite{Artyukhin1985}, a time varying boundary heat flux is studied. Other examples can also be consulted in \cite{Ucinski2004}.

An interesting remark with this approach is that the probability distribution of the unknown parameter $p_{\,m}$ can be estimated from the distribution of the measurements of the field $u$ and from the sensitivity $\Theta_{\,m}\,$. The probability $\mathcal{P}$ of $u$ is given by:
\begin{align*}
F \, (\, \bar{u}\,) \egal \mathcal{P} \, \biggl\{\, u(\,x \,, t \,, p_{\,m} \,) \ \leqslant \ \bar{u} \, \biggr\} \,.
\end{align*} 
Using the sensibility function $\Theta_{\,m}\,$, the probability can be approximated by:
\begin{align*}
  F \, (\, \bar{u} \,) & \simeq \ \mathcal{P} \, \biggl\{\, u(\,x \,, t \,, p_{\,m}^{\,\circ} \,) \plus \Theta_{\,m} \cdot  \bigl(\, p_{\,m} \moins p_{\,m}^{\,\circ} \,\bigr) \ \leqslant \ \bar{u} \, \biggr\} \,, \\
\end{align*} 
Assuming $\Theta_{\,m} \ > \ 0\,$, we get: 
\begin{align*}
  F \, (\, \bar{u} \,) & \egal  \mathcal{P} \, \biggl\{\, p_{\,m} \ \leqslant \  p_{\,m}^{\,\circ} \plus \frac{\bar{u} \moins u(\,x \,, t \,, p_{\,m}^{\,\circ} \,) }{\Theta_{\,m}}\, \biggr\} \,.
\end{align*}
Therefore, using a change o variable, the cumulative derivative function of the probability of the unknown parameter $p_{\,m}$ is estimated by:
\begin{align}\label{eq:probability_pm_1}
  F \, (\,\bar{p}_{\,m}\,) & \egal \mathcal{P} \, \biggl\{\, p_{\,m} \ \leqslant \ \bar{p}_{\,m} \, \biggr\} \,, \\
  & \egal F \, \biggl(\, u \plus \Theta_{\,m} \cdot \bigl(\, \bar{p}_{\,m} \moins p_{\,m}^{\,\circ} \,\bigr) \,\biggr) \,.
\end{align}
When $\Theta_{\,m} \ < \ 0\,$, the cumulative derivative function of the probability is given by:
\begin{align}\label{eq:probability_pm_2}
  F \, (\,\bar{p}_{\,m}\,) & \egal 1 \moins F \, \biggl(\, u \plus \Theta_{\,m} \cdot \bigl(\, \bar{p}_{\,m} \moins p_{\,m}^{\,\circ} \,\bigr) \,\biggr) \,.
\end{align}
It gives an \emph{local} approximation of the probability distribution of the unknown parameter $p_{\,m}\,$, at a reduced computation cost. Moreover, the approximation is reversible. Thus, if one has the distribution of the unknown parameter, it is possible to get the one of field $u\,$.


\section{Description of the experimental facility}
\label{sec:exp_facility}

\begin{figure}
\begin{center}
\subfigure[][\label{fig:schema_RH_box}]{\includegraphics[width=0.99\textwidth]{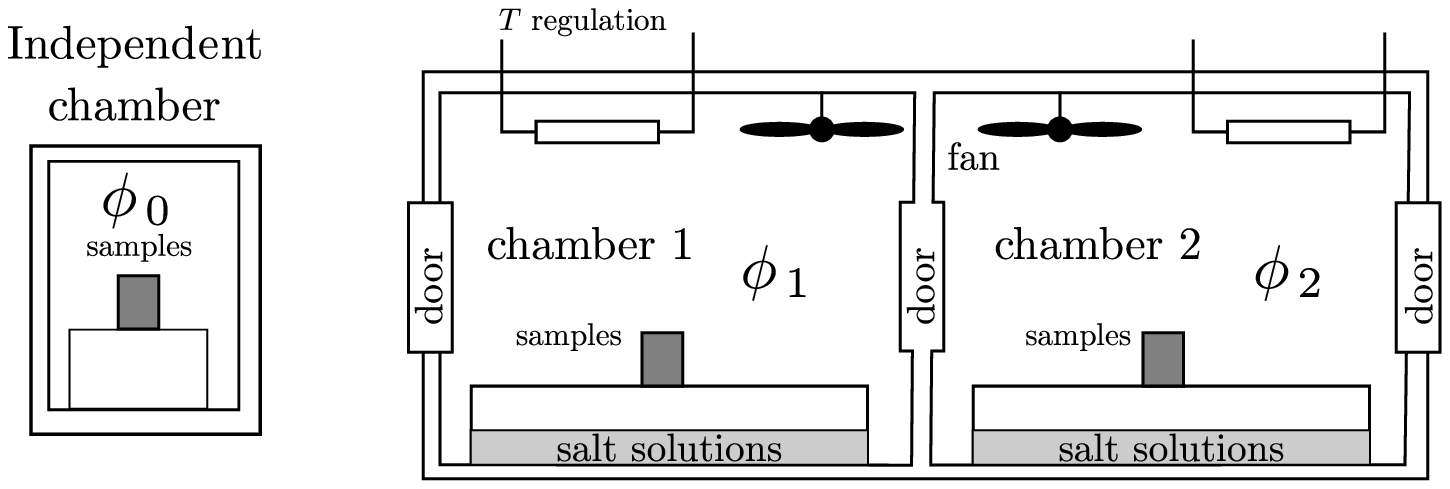}}
\subfigure[][\label{fig:photo_rh_box}]{\includegraphics[width=0.99\textwidth]{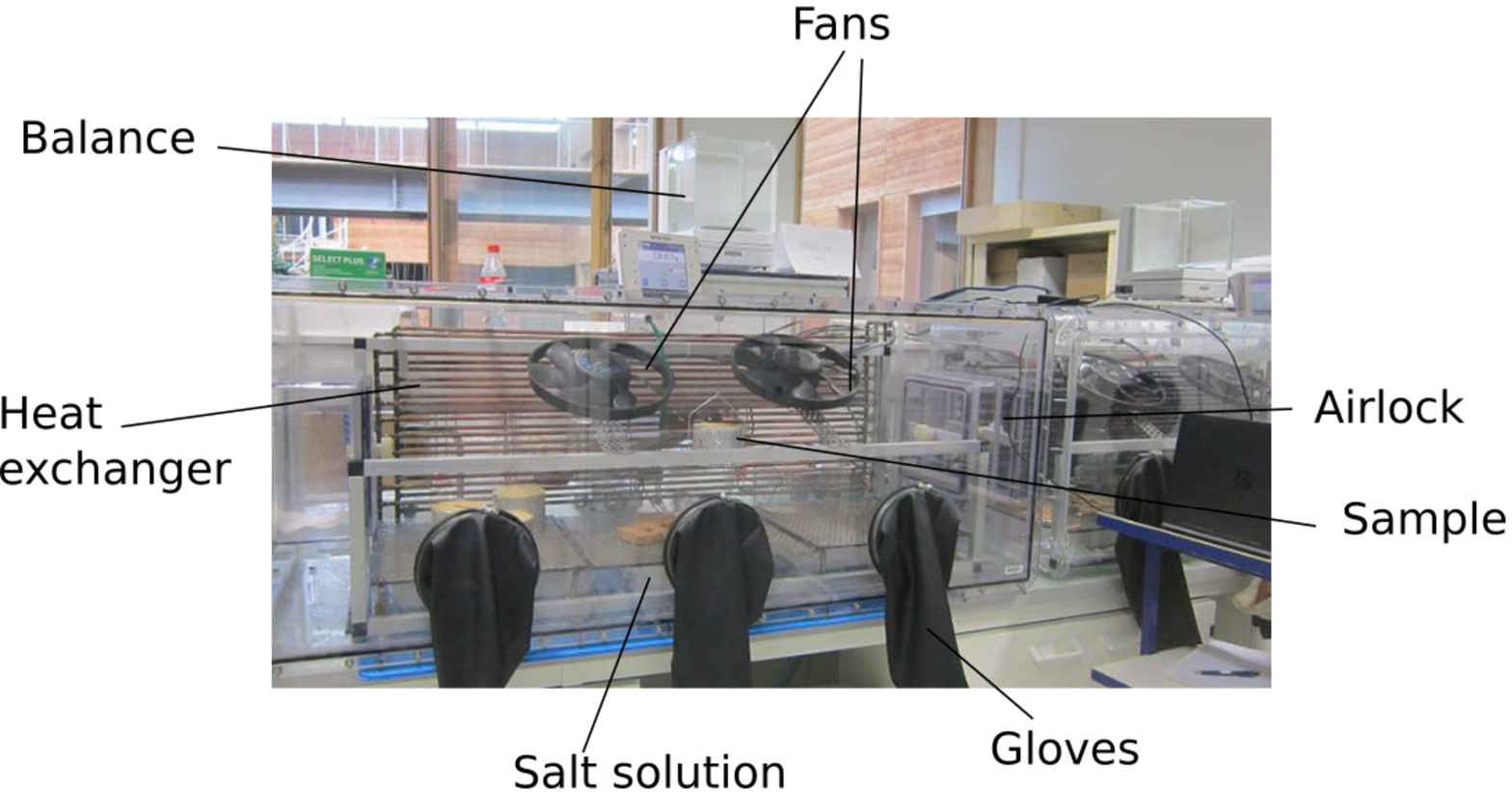}}
\caption{\small\em Illustration of the \emph{RH-box} experimental facility.}
\label{fig:RH_box}
\end{center}
\end{figure}

The test facility used to carry out the experiment, henceforth mentioned as the \emph{RH-Box}, is illustrated in Figure~\ref{fig:schema_RH_box} and a picture is given in Figure~\ref{fig:photo_rh_box}. It is composed of two connected climatic chambers. The temperature of each chamber is controlled independently with a thermostatically-controlled water bath allowing water to circulate in a heat exchanger. Relative humidity is controlled using $\mathrm{MgCl}_{\,2}$ and $\mathrm{NaCl}$ saturated salt solutions. Chambers $1$ and $2$ are fixed to $\phi_{\,1} \egal 33 \ \unit{\%}$ and $\phi_{\,2} \egal 75 \ \unit{\%}\,$, respectively. Another climatic chamber is also available. It is independent from the others as there is no doors linking the chambers. It is regulated at $\phi_{\,0} \egal 10 \ \unit{\%}$ to initially condition materials.

The samples are placed above the saline solution on grids. To avoid any perturbation of the chamber conditions, the various samples are handled using gloves. Two door locks, at each side, allow the operator to insert or remove samples to minimize system disturbances. They enable easy and instantaneous change in humidity boundary conditions for the samples while passing from one chamber to another. Two fans, operating continuously, allow to homogenize the inside air when the saline solutions are settled. Interested readers may consult a complementary discussion in \cite{Busser2016}.

\begin{figure}
\centering
\subfigure[][\label{fig:sensors}]{\includegraphics[width=0.35\textwidth]{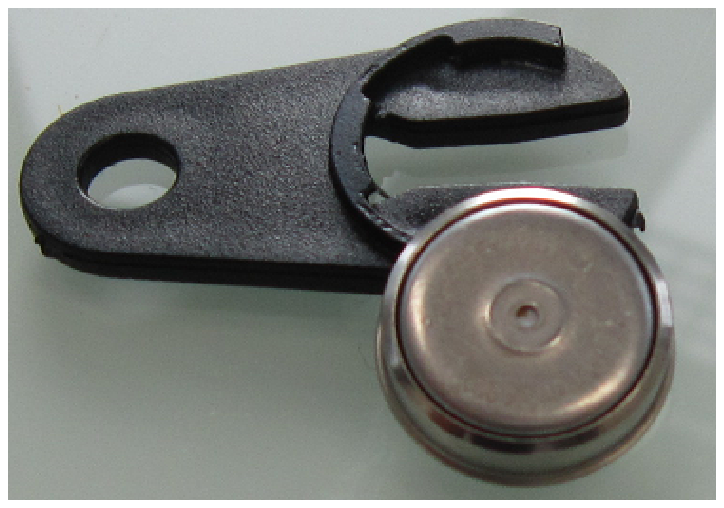}} 
\subfigure[][\label{fig:dessin_samples}]{\includegraphics[width=0.65\textwidth]{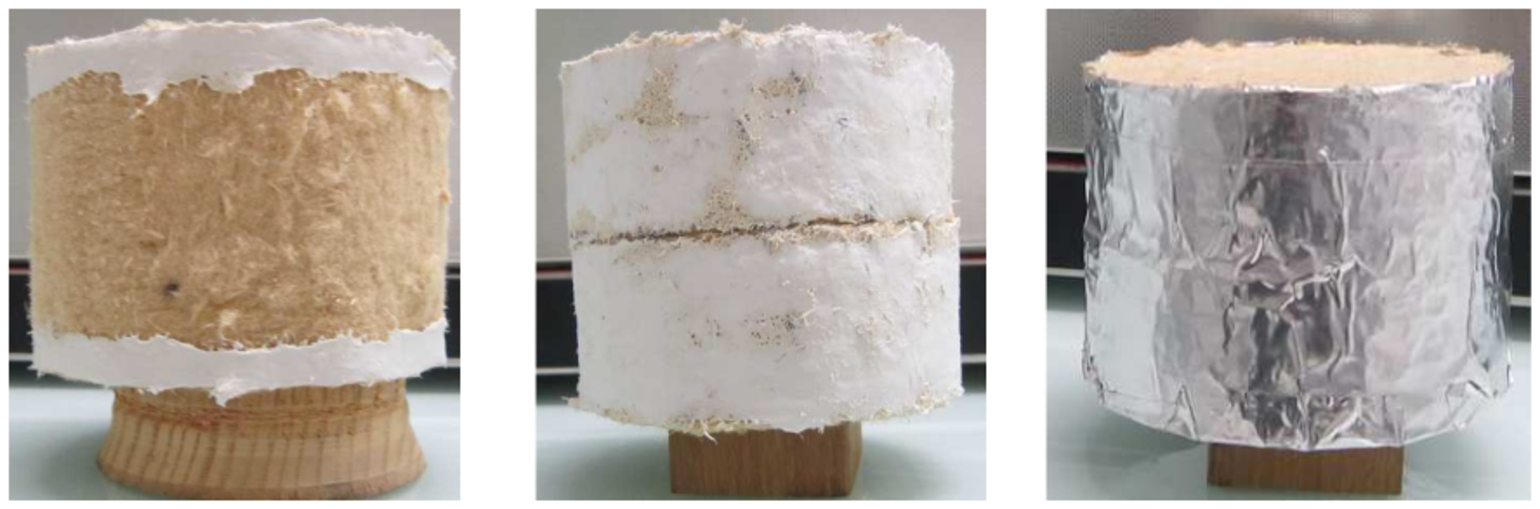}}
\caption{\small\em Sensors of relative humidity and temperature (a) and wood fibre samples (b) with white acrylic seal and with aluminium tape.}
\label{fig:samples_sensors}
\end{figure}

The temperature and relative humidity fields are measured inside the samples with wireless sensors from the \texttt{HygroPuce} range (\texttt{Waranet} industry). The accuracy is $\pm \ 2 \ \unit{\%}$ for the relative humidity and $\pm \ 0.5 \ \unit{^\circ C}$ for the temperature and the dimensions are $0.6 \ \unit{cm}$ thickness and $1.6 \ \unit{cm}$ diameter, as illustrated in Figure~\ref{fig:sensors}. The response time of the sensor is $5$ to $10 \ \mathsf{min}\,$.

\begin{figure}
\begin{center}
\subfigure[][\label{fig:TL_ft}]{\includegraphics[width=.48\textwidth]{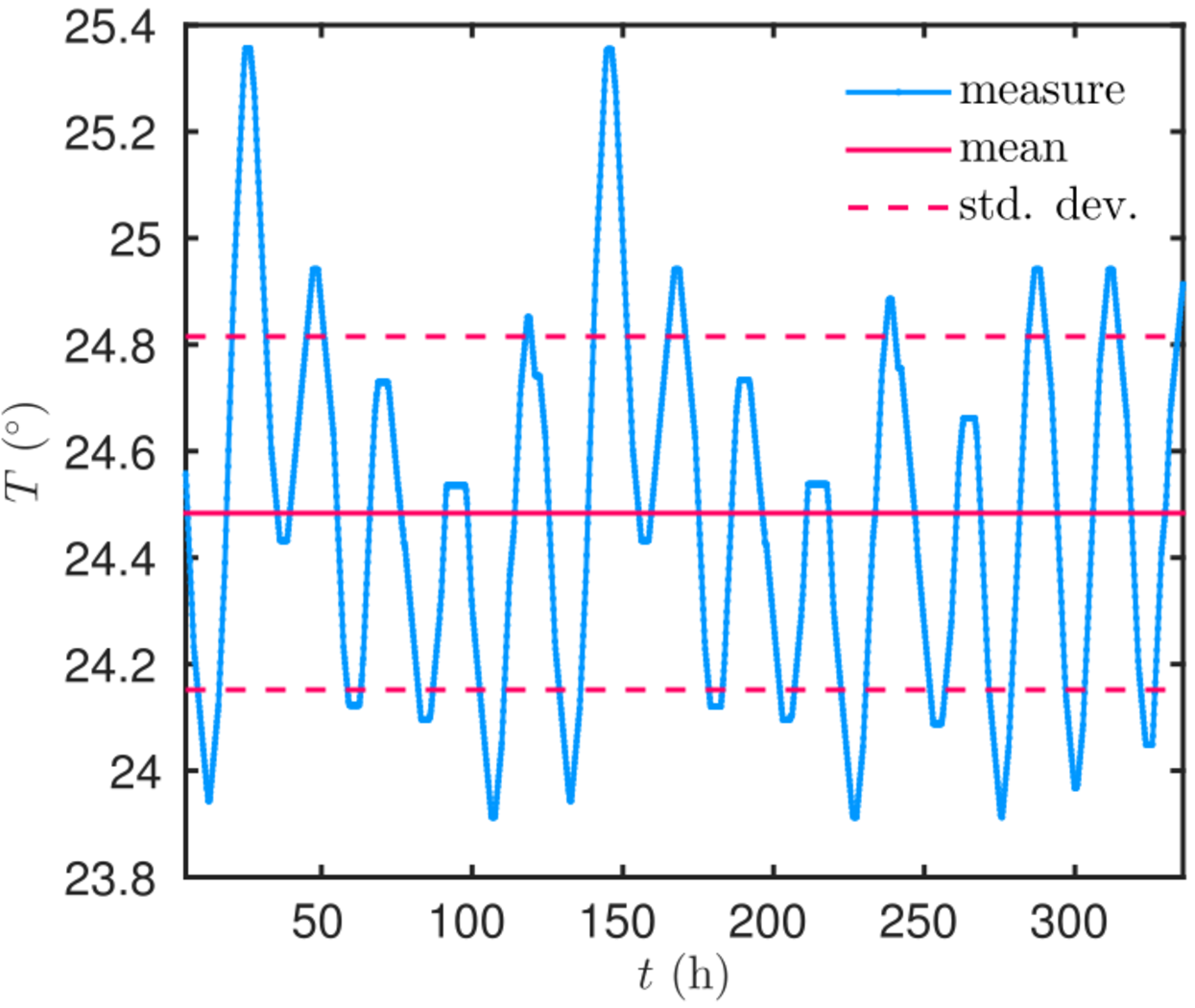}} \hspace{0.3cm}
\subfigure[][\label{fig:Tx3_ft}]{\includegraphics[width=.48\textwidth]{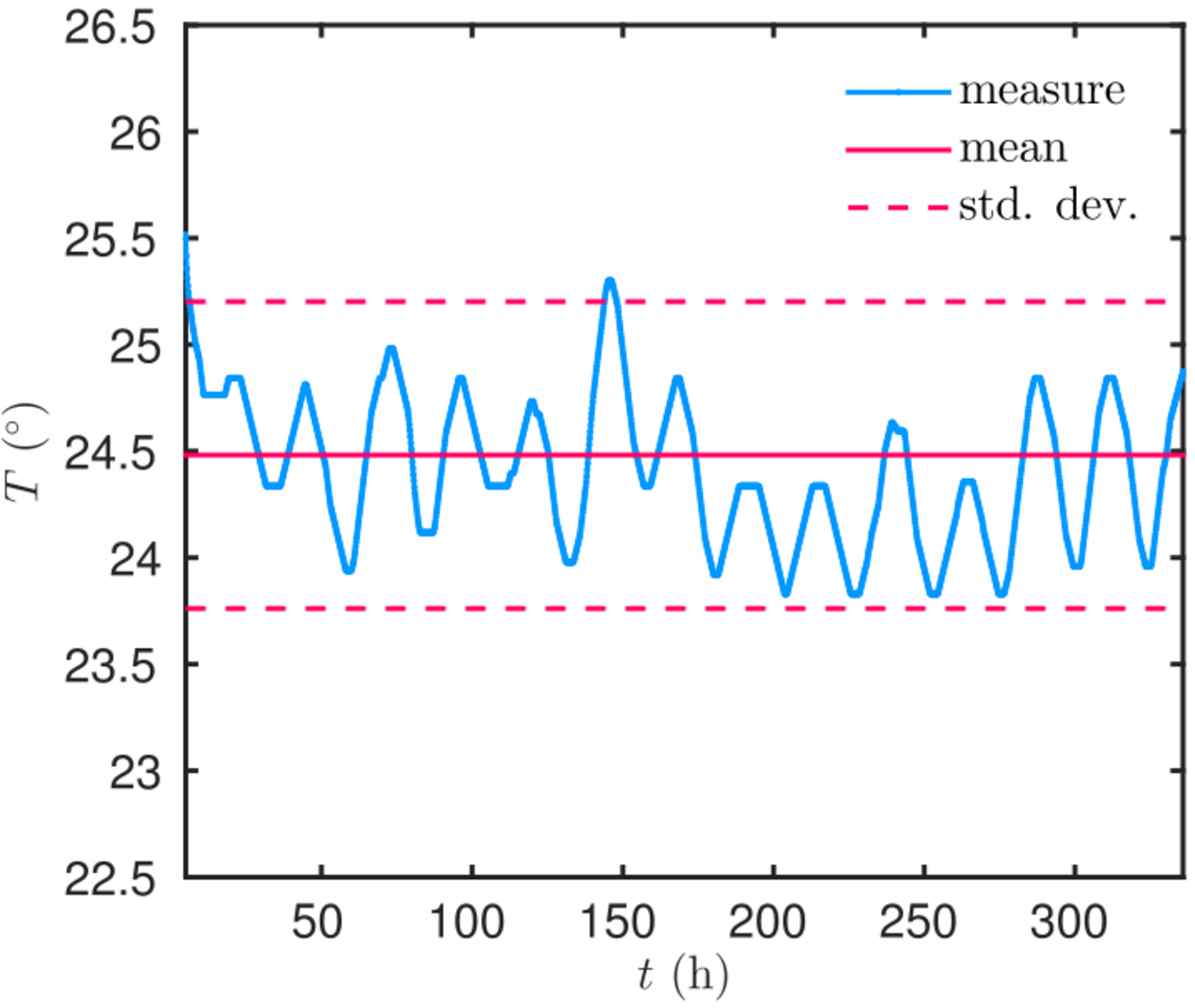}}
\caption{\small\em Variation of the ambient temperature inside the \emph{RH-box} (a) and inside the material at $x \egal 2.6 \ \unit{cm}$ (b).}
\label{fig:Tmeas_ft}
\end{center}
\end{figure}

The material investigated is the wood fibre, which properties have been determined in \cite{Rafidiarison2015, Vololonirina2014} and are shown in Table~\ref{tab:mat_properties}. They constitute \emph{a priori} information on the unknown parameters $d_{\,0}$ and $d_{\,1}\,$. The temperature was assumed to have slight variations in the material, around $T \egal 24.5 \ \unit{^\circ C}$. This hypothesis is assumed considering the slight variation of temperature observed in the ambient facility and in the material, as illustrated in Figures~\ref{fig:TL_ft} and \ref{fig:Tx3_ft}. In \cite{Berger2017a}, the mass average velocity in the material was estimated around $\mathsf{v} \egal 5 \ \unit{mm/s}\,$, where the facility imposed an air flow at one face of the material. In the current study, no air movement are imposed. Therefore a lower velocity value $\mathsf{v} \egal 0.1 \ \unit{mm/s}$ is used as \emph{a priori} value in the material. This hypothesis will be verified when estimating the parameters. Therefore, according to the definition, an \emph{a priori} estimation of the advection parameter is $a \egal 7.2 \e{-11} \ \unit{s/m}\,$. The samples are cylindrical, with a $10 \ \unit{cm}$ diameter and $8 \ \unit{cm} \, $ thickness. The size of the samples has been chosen in order to avoid borders effects and to minimize the perturbation of sensors inside the sample. Moreover, to ensure one-dimensional moisture transfer, the side and bottom surfaces of the samples are covered with aluminium tape and glued on a white acrylic seal, as illustrated in Figure~\ref{fig:dessin_samples}.

The sensors are placed within the wood fibre by cutting the samples. The exact location of the sensor will be investigated through the OED approach in the next Sections. The influence of cutting the sample and introducing the sensor inside the sample has been studied in \cite{Busser2016}. Several samples have been under investigation by cutting them and introducing $1\,$, $2$ or $3$ samples at different position. Sensors were inserted inside the sample and at the bottom surface of the sample. Measurement inside these samples has been compared to a reference measure. The latter was obtained from a sample with only a sensor at the bottom surface, settled without cutting the sample. Results revealed a good agreement between the measurements of all samples, concluding that the protocol may not disturb the moisture transfer.
The total uncertainty on the measurement of relative humidity can be evaluated considering the propagation of the uncertainty due to sensor measurement and due to their location. It should be noted that uncertainties due to contact resistance between the sensor and the material are difficult to estimate and therefore are not taken into account here. The total uncertainty on relative humidity is \cite{Taylor1997}:
\begin{align*}
  \Delta \phi \egal \phi \ \sqrt{ \Biggl(\, \pd{\phi}{x}\Biggr|_{\,x \egal X} \cdot \frac{\Delta X}{\phi}\,\Biggr)^{\,2} \plus \Biggl(\, \frac{\Delta \phi}{\phi} \,\Biggr)^{\,2} } \,,
\end{align*}
where $X$ is the location of the sensor and $\Delta X$ the uncertainty on its position. The term $\pd{\phi}{x}$ has been evaluated using a numerical computation of the direct problem for the different range of humidity investigated here and results show a maximum value of $0.04 \ \unit{m^{\,-1}}\,$. Considering measuring a relative humidity $\phi \egal 10 \ \unit{\%}\,$, the sensor uncertainty $\Delta \phi \egal 2 \ \unit{\%}\,$, and the sensor location uncertainty $\Delta X \egal 0.1 \ \unit{mm}\,$, the total uncertainty on the measurement is of the order of $\Delta \phi \egal 2 \ \unit{\%}\,$. The contribution of the uncertainty of the sensor location is $100$ times lower than the one of the sensor measurement. These results confirm the ones presented in \cite{Busser2016}, showing that the hygrothermal behaviour is not impacted by the experimental design. Moreover, the reproducibility of the experiments has also been verified, testing seven samples at the same time when carrying out the experiments.

\begin{figure}
\centering
\includegraphics[width=0.55\textwidth]{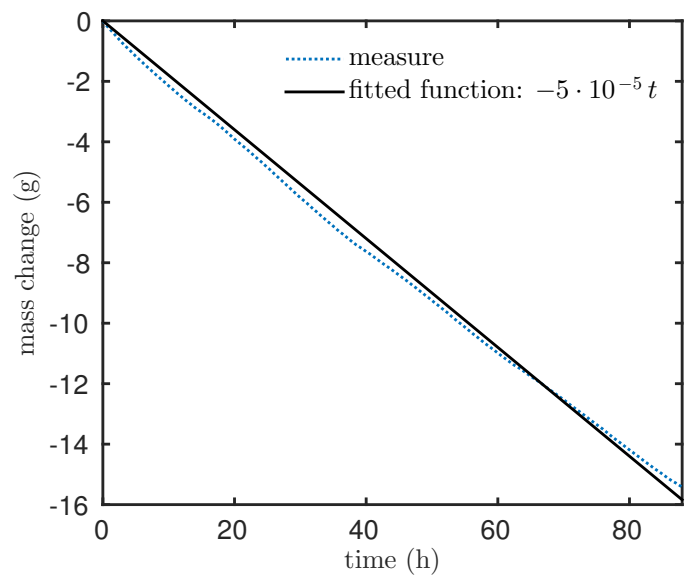}
\caption{\small\em Time variation of the mass changes for the estimation of the convective mass transfer coefficient $h\,$.}
\label{fig:mass_change}
\end{figure}

As mentioned in Section~\ref{sec:moisture_equations}, an null flux condition is assumed at the bottom of the samples ($x \egal L$). At $x \egal 0\,$, a \textsc{Robin} type boundary condition Eq.~\eqref{eq:HAM_BC_L} is considered. The surface is in contact with the ambient air of the chamber knowing that the experimental facility enables to perform one or multiple steps of relative humidity: 
\begin{align*}
  & \phiL \egal \left. \begin{cases}
    \phi^{\,i}  \,, & t \ < \ 0 \,,\\
    \phi^{\,i}  \plus  \dPhi \,, & t  \ \geqslant \ 0  \,,
  \end{cases} \right.
\end{align*}
The convective mass transfer coefficient $h$ has been estimated by preliminary measurements adapted from \cite{Kwiatkowski2009}. A reservoir with liquid water is placed in the \emph{RH-box} and the mass change was recorded during $\, 3.5$ days. Figure~\ref{fig:mass_change}  gives the time variations of water weight in the reservoir placed in the climatic chamber with $\phi \egal 33 \ \unit{\%}\,$. The convective mass transfer coefficient has been estimated using following equation:
\begin{align*}
  h \egal \frac{g}{ \Ps (\, T_{\,w} \,) \moins \Ps (\, \TL \,) \, \phiL \,} \,,
\end{align*}
where $\Ps (\, T_{\,w} \,)$ is the saturated vapour pressure for the temperature $T_{\,w}$ of the liquid water. The moisture flux $g$ is computed using the slope of the mass variation of water according to time. Results give a convective mass transfer coefficient around $h \egal 9 \e{-9} \ \unit{s/m}\,$. This value provides an acceptable estimation of the coefficient, given the values used in literature \cite{Rouchier2017, Olek2016}. It will be used as a first approximation for the hygroscopic material.

Finally, the experimental facility is used to submit the samples to a single or multiple steps of relative humidity. In the first case, the material can be initially conditioned at $\phi \egal 10 \ \unit{\%}\,$, $\phi \egal 33 \ \unit{\%}$ or $\phi \egal 75 \ \unit{\%}\,$. In the case of multiple steps, as the climatic chamber $\phi \egal 10 \ \unit{\%}$ is not linked with the others, the samples will always be conditioned using this chamber. Therefore, the multiple steps can only be operated between climatic chamber $1$ and $2\,$, where the perturbations, due to the mouvement of the samples, are reduced. For single step experiments, it is not possible to perform a step from $33 \ \unit{\%}$ or $75 \ \unit{\%}$ to the climatic chamber at $10 \ \unit{\%}$ for the same practical reasons. The pre-conditioning is reached when the mass change of the samples is lower than $0.05 \ \unit{\%}\,$. A total of $20$ designs are possible as reported in Tables~\ref{tab:1step_Design_1parameter} and \ref{tab:Mstep_Design_1parameter}.

\begin{table}
\centering
\setlength{\extrarowheight}{.3em}
\begin{tabular}[l]{@{} lc}
\hline
\hline
Sorption curve $\biggl[\, \unitfrac{kg\,}{\,m^{\,3}} \,\biggr]$ & $\f(\phi) \egal 40 \ \phi^{\,3} - 49 \ \phi^{\,2} + 27.02 \ \phi$ \\
Moisture permeability $\bigl[\, \unit{s} \,\bigr]$ & $d \egal 2.33 \e{-11} \plus 5.68 \e{-11} \ \phi$\\
Advection coefficient $\bigl[\, \unit{s/m} \,\bigr]$ & $a \egal 7.2 \e{-11}$ \\
\hline
\hline
\end{tabular}
\bigskip
\caption{\small\em \emph{A priori} material properties of wood fibre from \cite{Rafidiarison2015, Berger2017a}.}
\label{tab:mat_properties}
\end{table}


\section{Determination of the OED}
\label{sec:OED}

\subsection{Estimation of one parameter}

\subsubsection{Single step of relative humidity}

The aim is to use the measurements to estimate one of the parameters of the material moisture permeability and advection coefficients, $d_{\,0}\,$, $d_{\,1}$ or $a\,$, with a single step in relative humidity. Considering the facility described in Section~\ref{sec:exp_facility}, four designs are possible for a single step of humidity, and are reported in Table~\ref{tab:1step_Design_1parameter}. The issue is to determine the optimal design, in terms of location of the sensors in the material for each of the four designs. Then, by comparing the magnitude of the criterion $\Psi\,$, one can define which step of relative humidity is the best design. For this purpose, the duration of experiments is set to $\tau \egal 200 \ \unit{h}\,$. These wood material properties reported in Table~\ref{tab:mat_properties} are used as prior informations to compute the sensitivity coefficients.

\begin{figure}
\begin{center}
\subfigure[][\label{fig:1step_Psi_1param}]{\includegraphics[width=0.47\textwidth]{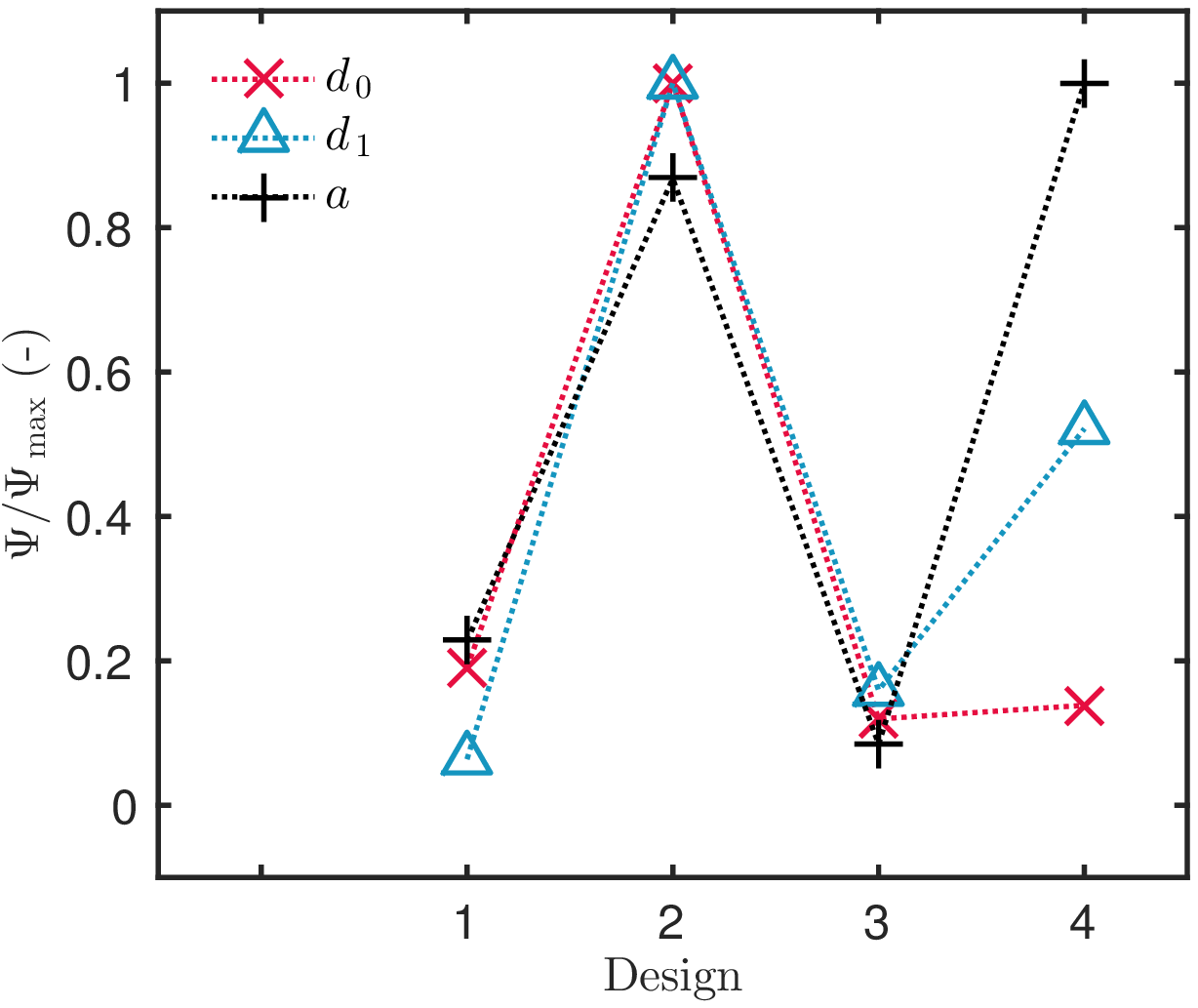}}
\subfigure[][\label{fig:1step_Psi_fx_1param}]{\includegraphics[width=0.49\textwidth]{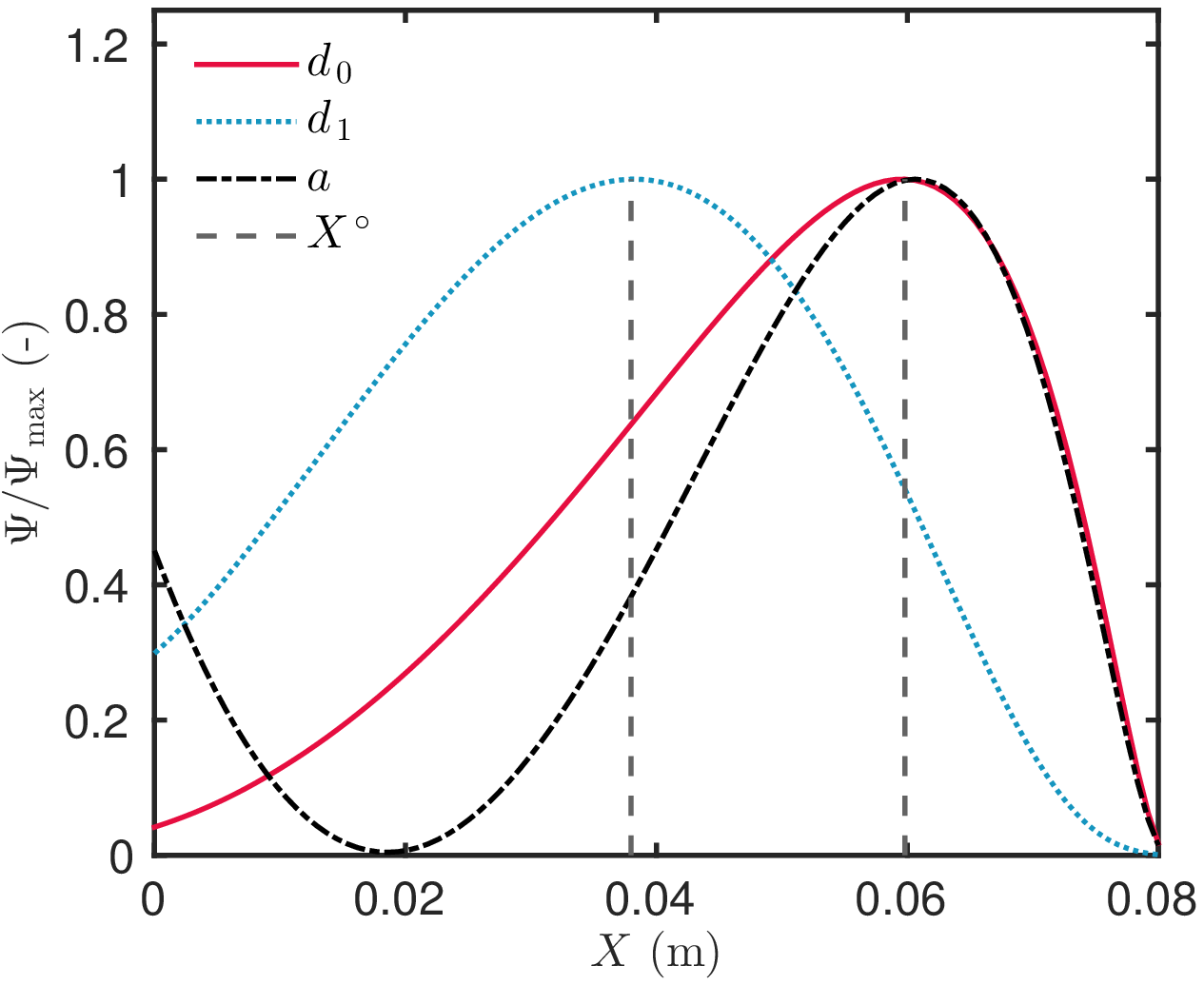}}
\caption{\small\em Variation of the criterion $\Psi$ for the four possible designs (a) and as a function of the sensor position $X$ for design 2 (b).}
\label{fig:1step_Psi_design}
\end{center}
\end{figure}

\begin{figure}
\begin{center}
\subfigure[][\label{fig:1step_S_OED}]{\includegraphics[width=.45\textwidth]{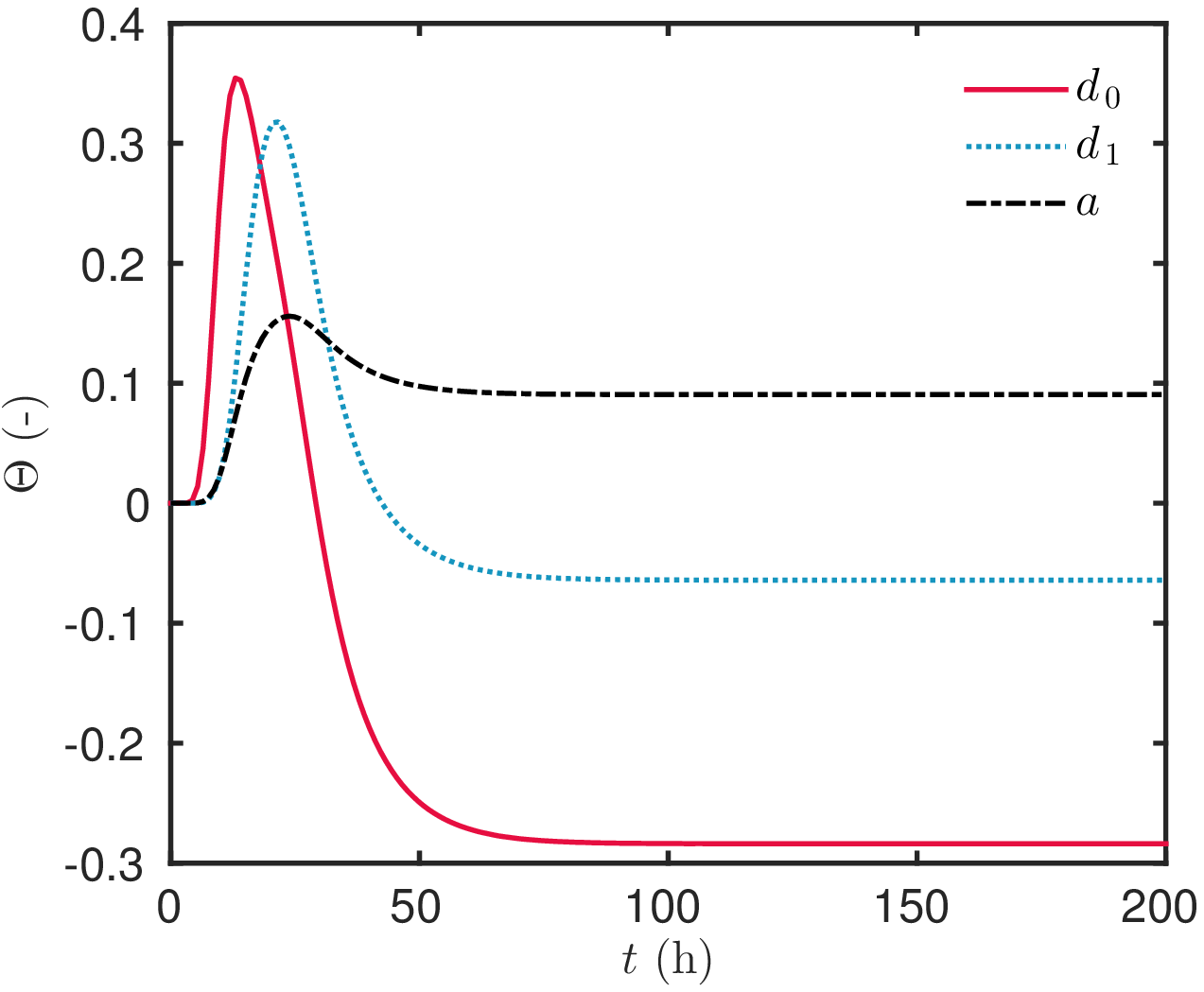}}
\subfigure[][\label{fig:1step_S_anti_OED}]{\includegraphics[width=.45\textwidth]{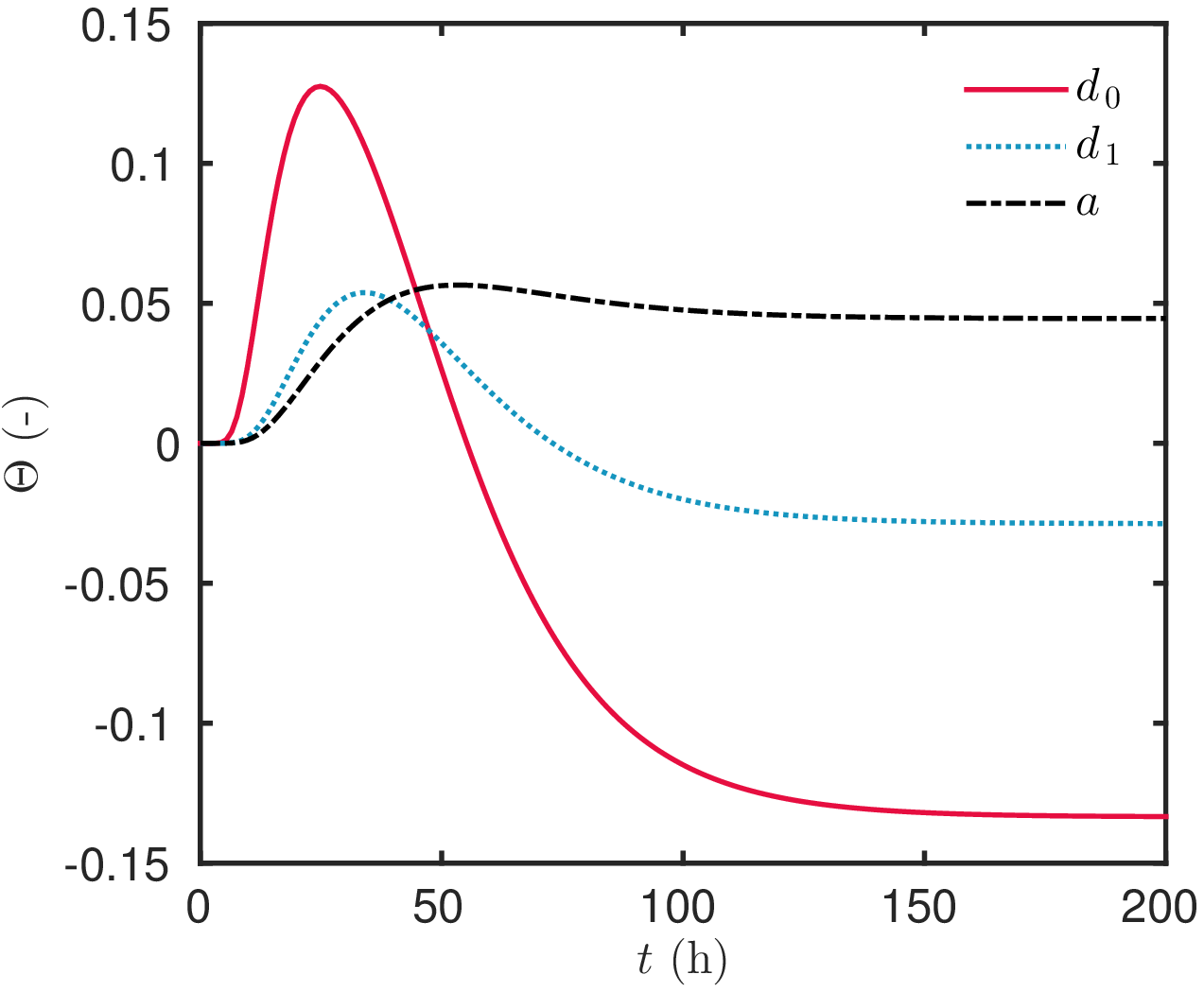}}
\caption{\small\em Sensitivity coefficients $\Theta$ for parameters $d_{\,0}\,$, $d_{\,1}$ and $a$ for the designs $2$ (OED) (a) and $1$ (c) ($X \egal X^{\,\circ}$).}
\label{fig:1step_S_OED_anti_OED}
\end{center}
\end{figure}

\begin{figure}
\begin{center}
\includegraphics[width=0.55\textwidth]{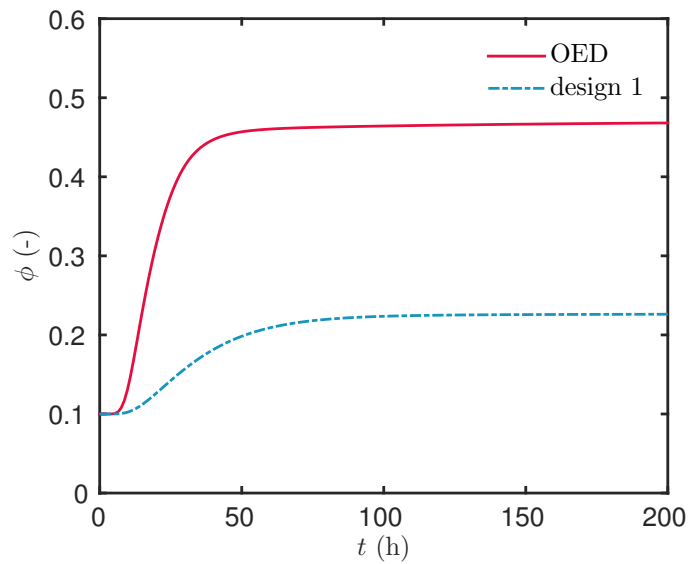}
\caption{Variation of the relative humidity $\phi$ for the designs $2$ (OED) (a) and $1$ (c) ($X \egal X^{\,\circ}$).}
\label{fig:1step_v_OED_antiOED}
\end{center}
\end{figure}

The $D-$criterion $\Psi$ has been computed for each of the four designs and results are given in Table~\ref{tab:1step_Design_1parameter}. As illustrated in Figure~\ref{fig:1step_Psi_1param}, the design $2\,$, which has the higher step of relative humidity, has the highest criterion $\Psi$ for both parameters $d_{\,0}$ and $d_{\,1}\,$. This design would provide the highest accuracy to estimate one of these two parameters. For the advection coefficient $a\,$, the design $4$ has the criterion. Nevertheless, it can be noted that the design $2$ reaches almost $90\%$ of the maximal criterion. Therefore, both designs $2$ and $4$ could be used to estimate the parameter $a$ with the best accuracy. On the opposite, the designs $1$ and $3$ are the worst protocols. The results are confirmed by the variation of the sensitivity coefficients given in Figures~\ref{fig:1step_S_OED} and \ref{fig:1step_S_anti_OED}. The function $\Theta$ measures the sensitivity of the estimated field with respect to change in the parameter ($d_{\,0} \,$, $d_{\,1}\,$, or $a$). The function have higher magnitudes of variations for the design $2$ than for design $1\,$. The field is thus more sensible to changes in those parameters for the design $2\,$. Indeed, the relative humidity in the material have sharper variations in the material for the OED than for the worst design, as illustrated in Figure~\ref{fig:1step_v_OED_antiOED}. As noticed in Eq.~\eqref{sec1_eq:fisher_matrix}, the criterion $\Psi$ is computed using the sensitivity function $\Theta\,$, which explains why the OED corresponds to the design $2\,$.

The variation of the criterion $\Psi$ with the position $X$ of the sensor is given in Figure~\ref{fig:1step_Psi_fx_1param} for the design $2\,$. The criterion reaches its maximum for the sensor located between $X^{\,\circ} \egal 4 \ \unit{cm}$ and $X^{\,\circ} \egal 6 \ \unit{cm}$. This optimal position does not vary much for other designs for parameters $d_{\,0}$ and $d_{\,1}$), as reported in Table~\ref{tab:1step_Design_1parameter}. For parameter $a\,$, the optimal position varies between $3 \ \unit{cm}$ and $6 \ \unit{cm}\,$.

Summarizing, for the identification of parameter $d_{\,0}$ or $d_{\,1}\,$, a single step of humidity from $10 \ \unit{\%}$ to $75 \ \unit{\%}\,$, with a sensor located at $X^{\,\circ} \egal 6 \ \unit{cm}$ and $X^{\,\circ} \egal 4 \ \unit{cm}\,$, respectively, ensures the highest accuracy for solving the parameter estimation problem. Concerning the advection coefficient $a\,$, a step from $75 \ \unit{\%}$ to $33 \ \unit{\%}\,$, with a sensor located at $X^{\,\circ} \egal 3 \ \unit{cm}$ is the optimal design. The same design as the one for parameter $d_{\,1}$ can also be used to estimate this parameter.

\begin{table}
\centering
\setlength{\extrarowheight}{.3em}
\begin{tabular}[l]{@{} lcccccccc}
\hline
\hline
\multirow{2}{*}{Design} 
& Relative humidity 
& \multicolumn{2}{c}{$d_{\,0}$} 
& \multicolumn{2}{c}{$d_{\,1}$} 
& \multicolumn{2}{c}{$a$} \\
&  step $ \phi \ (\,\unit{\%} \,)$ 
& $\Psi_{\,\mathrm{max}}$ & $X^{\,\circ} \ (\,\unit{m}\,)$ 
& $\Psi_{\,\mathrm{max}}$ & $X^{\,\circ} \ (\,\unit{m}\,)$ 
& $\Psi_{\,\mathrm{max}}$ & $X^{\,\circ} \ (\,\unit{m}\,)$ \\
\hline
1 & 10-33 & 0.58 & 0.050 & 0.09 & 0.045 & 0.19 & 0.054 \\
2 & 10-75 & \textbf{3.08} & \textbf{0.058} & \textbf{1.48} & \textbf{0.039} & 0.73 & 0.059\\
3 & 33-75 & 0.36 & 0.050 & 0.23 & 0.043 & 0.07 & 0.058\\
4 & 75-33 & 0.42 & 0.040 & 0.77 & 0.040 & \textbf{0.84} & \textbf{0.031} \\
\hline
\hline
\end{tabular}
\bigskip
\caption{\small\em Determining the OED (highlighted in bold) for the identification of one parameter, using an experiment with a single step of relative humidity.}
\label{tab:1step_Design_1parameter}
\end{table}

\subsubsection{Multiple steps of relative humidity}

The OED has been obtained for the estimation of one parameter, considering an experiment with multiple steps of relative humidity. According to the description of the facility and to Table~\ref{tab:Mstep_Design_1parameter}, $16$ designs are considered. The design with a primary step to $\phi \egal 75 \unit{\%}$ has the highest criterion $\psi$ and thus is the OED. Nevertheless, it can be noted that this experiment require a long-time measurements (around $500 \ \unit{h}$). Therefore, the choice of an experiment with a single step relative humidity is preferable for the estimation of one parameter. In terms of sensor location, the OED is reached for $X^{\,\circ} \egal 5.5 \ \unit{cm}\,$, $X^{\,\circ} \egal 3.8 \ \unit{cm}$  and $X^{\,\circ} \egal 6 \ \unit{cm}$ for the estimation of $d_{\,0}\,$, $d_{\,1}$ and $a$ respectively.

\begin{table}
\centering
\setlength{\extrarowheight}{.3em}
\begin{tabular}[l]{@{} lccccccc}
\hline
\hline
\multirow{2}{*}{Design}  
& Duration of the  
& \multicolumn{2}{c}{$d_{\,0}$} 
& \multicolumn{2}{c}{$d_{\,1}$} 
& \multicolumn{2}{c}{$a$} \\
& steps $(\,\unit{days} \,)$
& $\Psi$ & $X^{\,\circ} \ (\,\unit{m}\,)$ 
& $\Psi$ & $X^{\,\circ} \ (\,\unit{m}\,)$ 
& $\Psi$ & $X^{\,\circ} \ (\,\unit{m}\,)$  \\
\hline
 \multicolumn{8}{c}{step $ \phi \ (\,\unit{\%} \,)$: 10-33-75-33} \\
1 & 1
& 0.49 & 0.055 
& 0.24 & 0.049
& 0.48 & 0.057  \\
2 & 2
& 0.83 & 0.056 
& 0.34 & 0.045
& 0.77 & 0.059  \\
3 & 3
& 1.13 & 0.058 
& 0.41 & 0.042 
& 1.03 & 0.059 \\
4 & 4
& 1.35 & 0.057 
& 0.45 & 0.042 
& 1.22 & 0.057 \\
5 & 5
& 1.48 & 0.056 
& 0.48 & 0.042 
& 1.33 & 0.060 \\
6 & 6
& 1.49 & 0.058 
& 0.48 & 0.042
& 1.34 & 0.060  \\
7 & 7 
& 1.49 & 0.058 
& 0.48 & 0.040
& 1.34 & 0.060 \\
8 & 8
& 1.49 & 0.058 
& 0.48 & 0.042 
& 1.36 & 0.060\\
\hline
\multicolumn{8}{c}{step $ \phi \ (\,\unit{\%} \,)$: 10-75-33-75} \\
9 & 1 
& 1.91 & 0.058 
& 1.05 & 0.042
& 1.65 & 0.061 \\
10 & 2
& 3.1 & 0.056 
& 1.82 & 0.039
& 4.3 & 0.061  \\
11 & 3 
& 5.3 & 0.056 
& 2.43 & 0.038
& 2.57 & 0.061 \\
12 & 4 
& 6.01 & 0.056 
& 2.92 & 0.038 
& 3.37 & 0.061 \\
13 & 5 
& 6.7 & 0.056 
& 3.28 & 0.038 
& 4.12 & 0.061\\
14 & 6 
& 6.9 & 0.056 
& 3.51 & 0.038 
& 1.58 & 0.061 \\
15 & 7
& 7.0 & 0.055 
& 3.56 & 0.038 
& 5.02 & 0.061 \\
\textbf{16} & \textbf{8} 
& \textbf{7.1} & \textbf{0.056}
& \textbf{3.8} & \textbf{0.038}
& \textbf{5.14} & \textbf{0.06} \\
\hline
\hline
\end{tabular}
\bigskip
\caption{\small\em Determining the OED (highlighted in bold) for the identification of one parameter, using an experiment with a multiple steps of relative humidity.}
\label{tab:Mstep_Design_1parameter}
\end{table}


\subsection{Estimation of several parameters}

\subsubsection{Single step relative humidity}

\begin{figure}
\begin{center}
\subfigure[][\label{fig:1step_Psi_2param}]{\includegraphics[width=.45\textwidth]{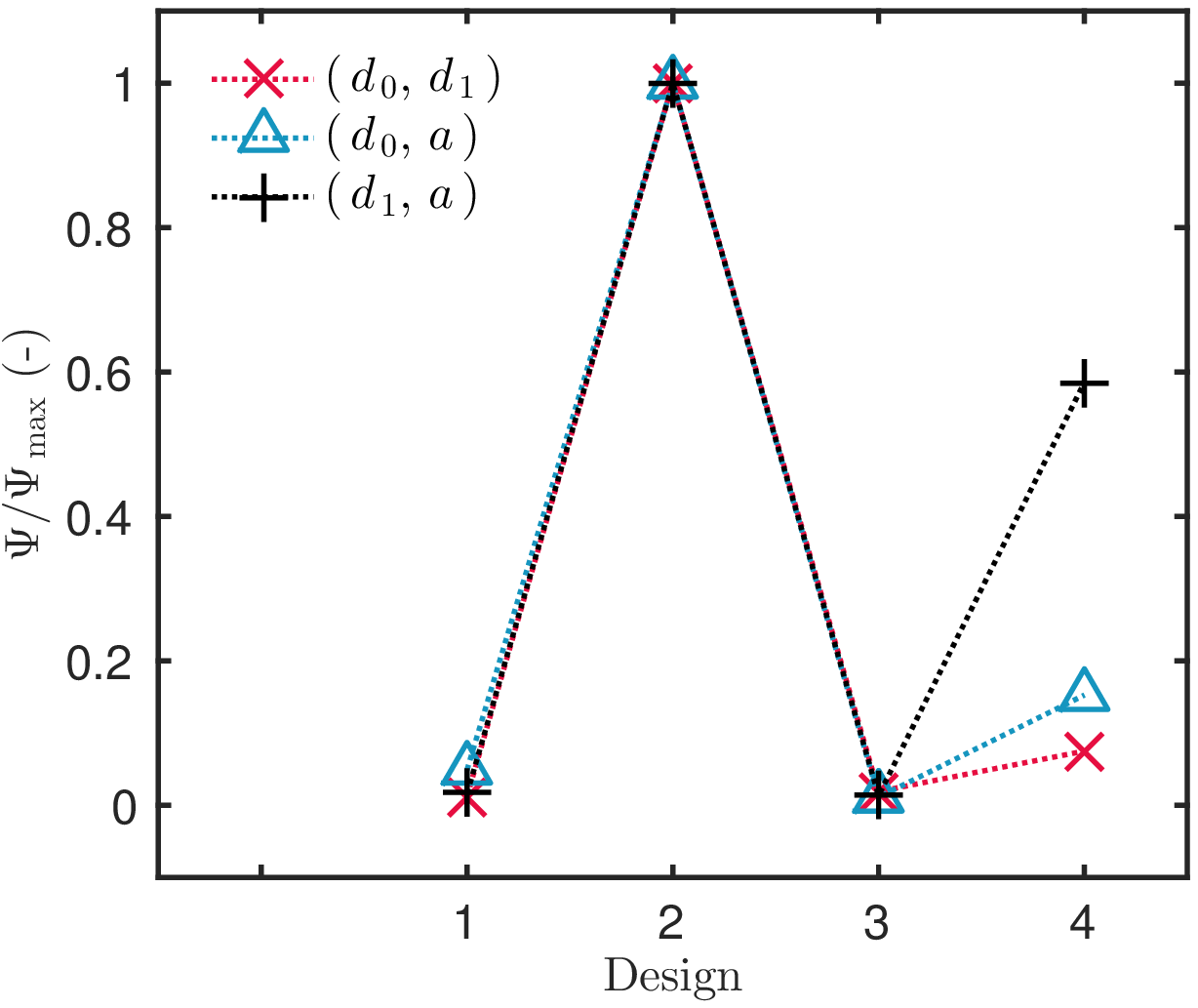}} \hspace{0.3cm}
\subfigure[][\label{fig:1step_Psi_fx_2param}]{\includegraphics[width=.45\textwidth]{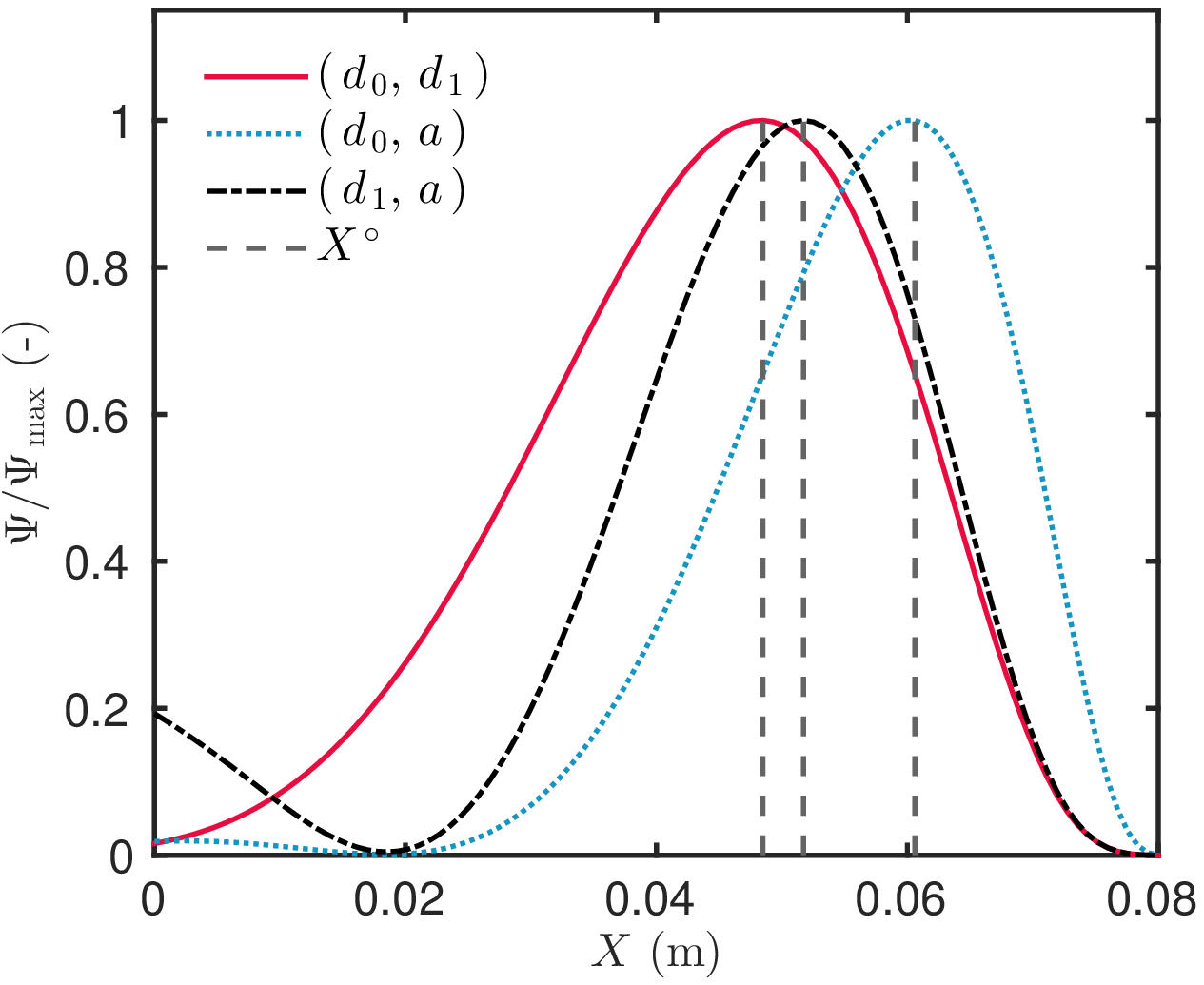}}
\caption{\small\em Variation of the criterion $\Psi$ for the four possible designs (a) and as a function of the sensor position $X$ for design $2$ (b).}
\label{fig_AN1:Psi_design_1}
\end{center}
\end{figure}

The issue is now to estimate two or three parameters among the material moisture permeability and advection coefficients, $d_{\,0}\,$, $d_{\,1}$ and $a\,$. This section focuses, therefore, on a single step of relative humidity. Table~\ref{tab:1step_Design_2parameters} gives an overview of the results. First of all, the sensitivities of parameters $d_{\,0}$ and $d_{\,1}$ have a strong correlation. This interaction can also be observed in Figure~\ref{fig:1step_S_OED_anti_OED}. Therefore, the identification of the couple $\bigl(\, d_{\,0},\, d_{\,1}\,\bigr)$ or even of three parameters might be a difficult task with these experiences. As shown in Figure~\ref{fig:1step_Psi_fx_2param}, the design $2$ appears to be the one providing the optimal conditions to estimate the couples of parameters $(\,d_{\,0} ,\, a \,)$ and $(\, d_{\,1} ,\, a \,)$. Figure~\ref{fig:1step_Psi_fx_2param} gives the optimal localisation of the sensors for the design $2\,$. The sensor should be located around $X^{\,\circ} \egal 6 \ \unit{cm}$ and $X^{\,\circ} \egal 5 \ \unit{cm}$  for the couples of parameters $(\,d_{\,0} ,\, a \,)$ and $(\, d_{\,1} ,\, a \,)\,$, respectively.

\begin{table}
\centering
\setlength{\extrarowheight}{.3em}
\begin{tabular}[l]{@{} lcccccccccccc}
\hline
\hline
\multirow{2}{*}{Design}  
& Relative humidity 
& \multicolumn{3}{c}{$(\,d_{\,0} ,\, d_{\,1} \,)$} 
& \multicolumn{3}{c}{$(\,d_{\,0} ,\, a \,)$} 
& \multicolumn{3}{c}{$(\, d_{\,1} ,\, a \,)$} \\
& step $ \phi \ (\,\unit{\%} \,)$ 
& $\Psi$ & $X^{\,\circ} \ (\,\unit{m}\,)$ & $\mathrm{Corr}$
& $\Psi$ & $X^{\,\circ} \ (\,\unit{m}\,)$ & $\mathrm{Corr}$
& $\Psi$ & $X^{\,\circ} \ (\,\unit{m}\,)$ & $\mathrm{Corr}$\\
\hline
1 & 10-33 
& 0.05 & 0.050 & 0.81
& 0.11 & 0.051 & -0.1
& 0.02 & 0.050 & 0.50\\
2 & 10-75 
& \textbf{4.40} & \textbf{0.048} & \textbf{0.85}
& \textbf{2.25} & \textbf{0.061} & \textbf{-0.16}
& \textbf{1.15} & \textbf{0.052} & \textbf{0.50}\\
3 & 33-75 
& 0.08 & 0.046 & 0.89
& 0.02 & 0.054 & -0.45
& 0.01 & 0.053 & 0.03\\
4 & 75-33 
& 0.33 & 0.040 & 0.83
& 0.35 & 0.030 & 0.045
& 0.61 & 0.041 & -0.46\\
\hline
\hline
\end{tabular}
\bigskip
\caption{\small\em Determining the OED (highlighted in bold) for the identification of several parameters, using an experiment with a single step of relative humidity.}
\label{tab:1step_Design_2parameters}
\end{table}


\subsubsection{Multiple steps of relative humidity}

\begin{figure}
\begin{center}
\subfigure[][\label{fig:Mstep_Psi_2param}]{\includegraphics[width=.48\textwidth]{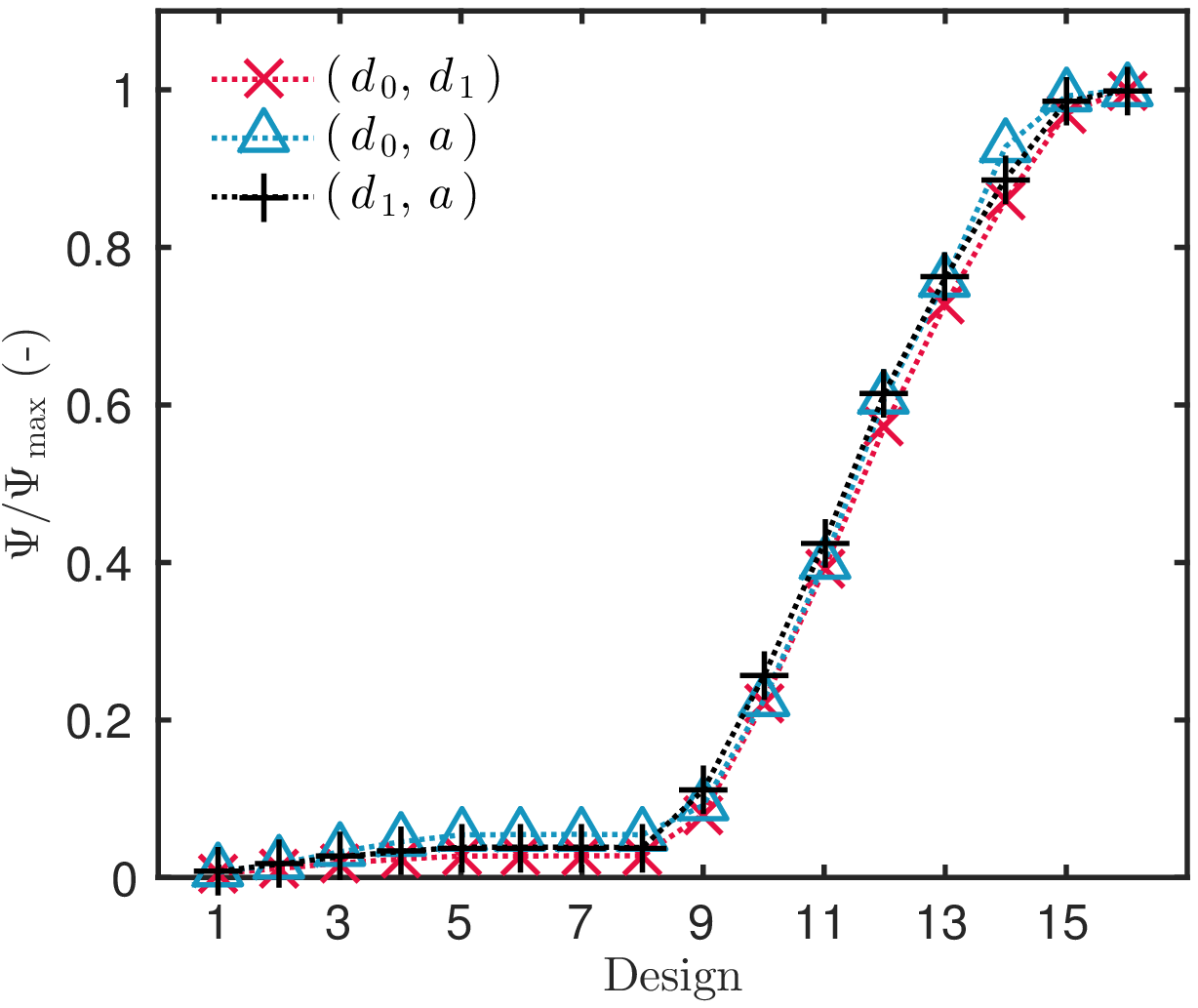}}
\subfigure[][\label{fig:Mstep_Psi_fx_2param}]{\includegraphics[width=.48\textwidth]{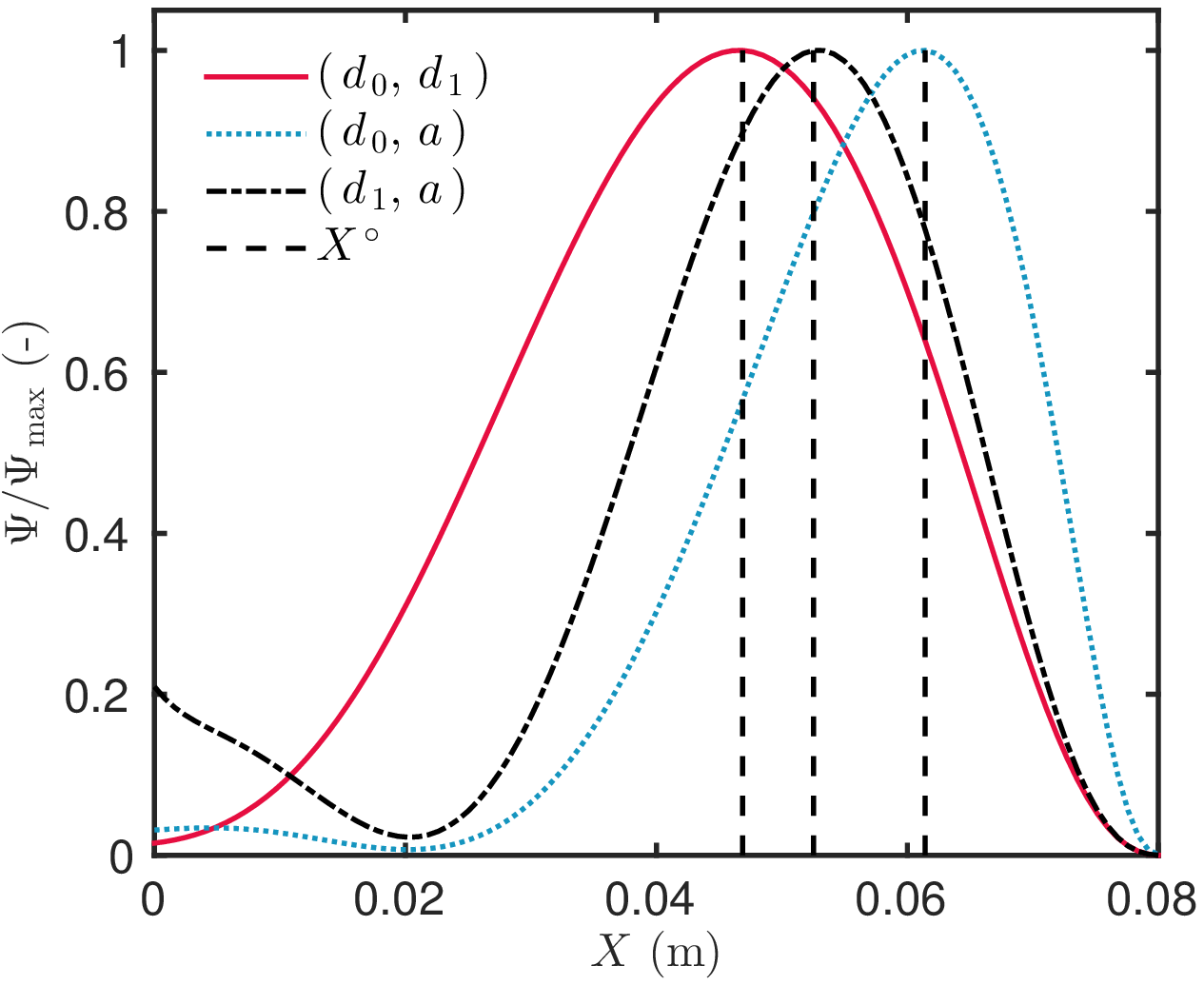}} \\
\caption{\small\em Variation of the criterion $\Psi$ for the four possible designs and for the identification of one (a) or two parameters (b).}
\label{fig:Mstep_Psi_design}
\end{center}
\end{figure}

\begin{figure}
\begin{center}
\subfigure[][\label{fig:Mstep_S_OED}]{\includegraphics[width=0.48\textwidth]{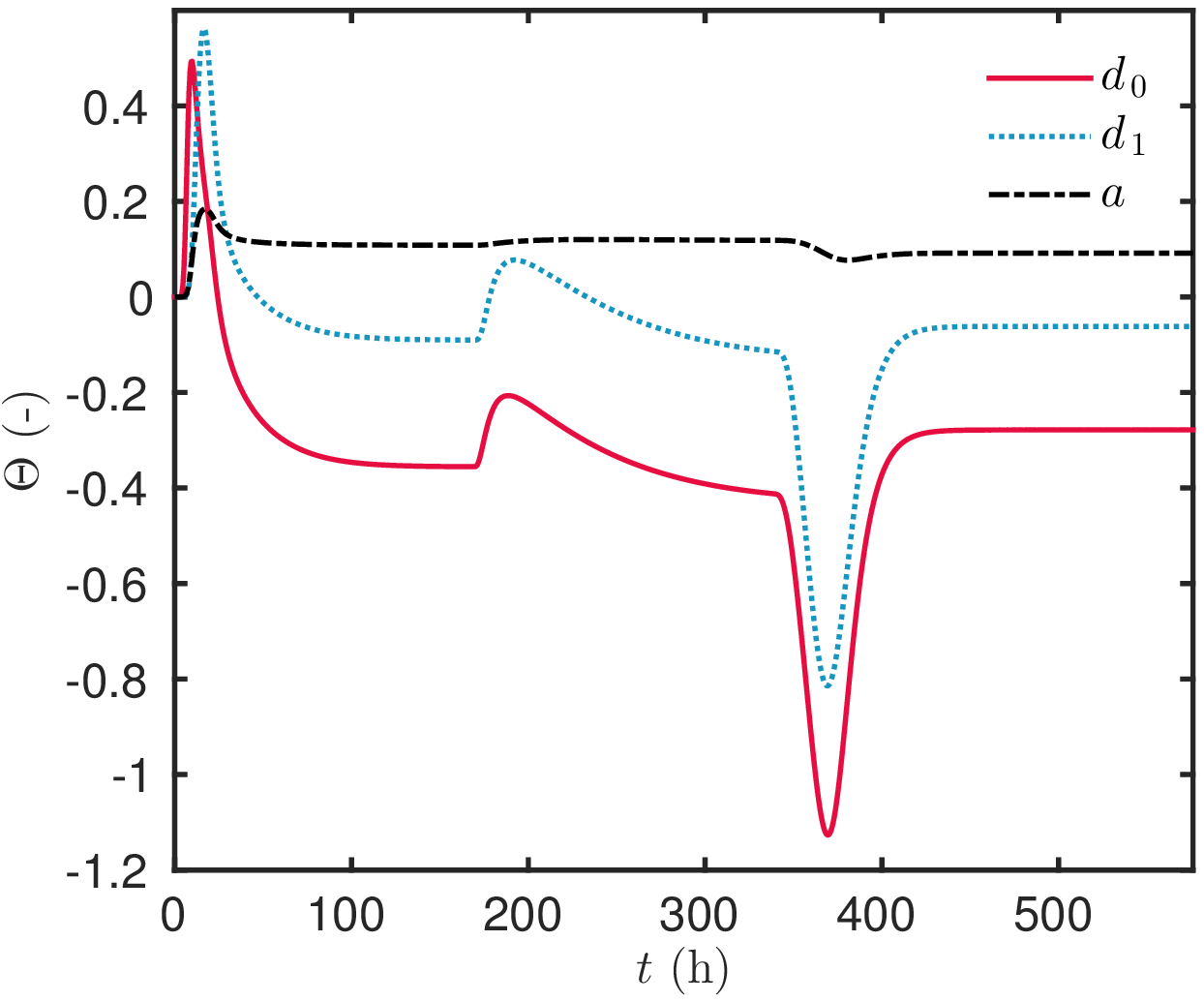}}
\subfigure[][\label{fig:Mstep_S_8}]{\includegraphics[width=0.48\textwidth]{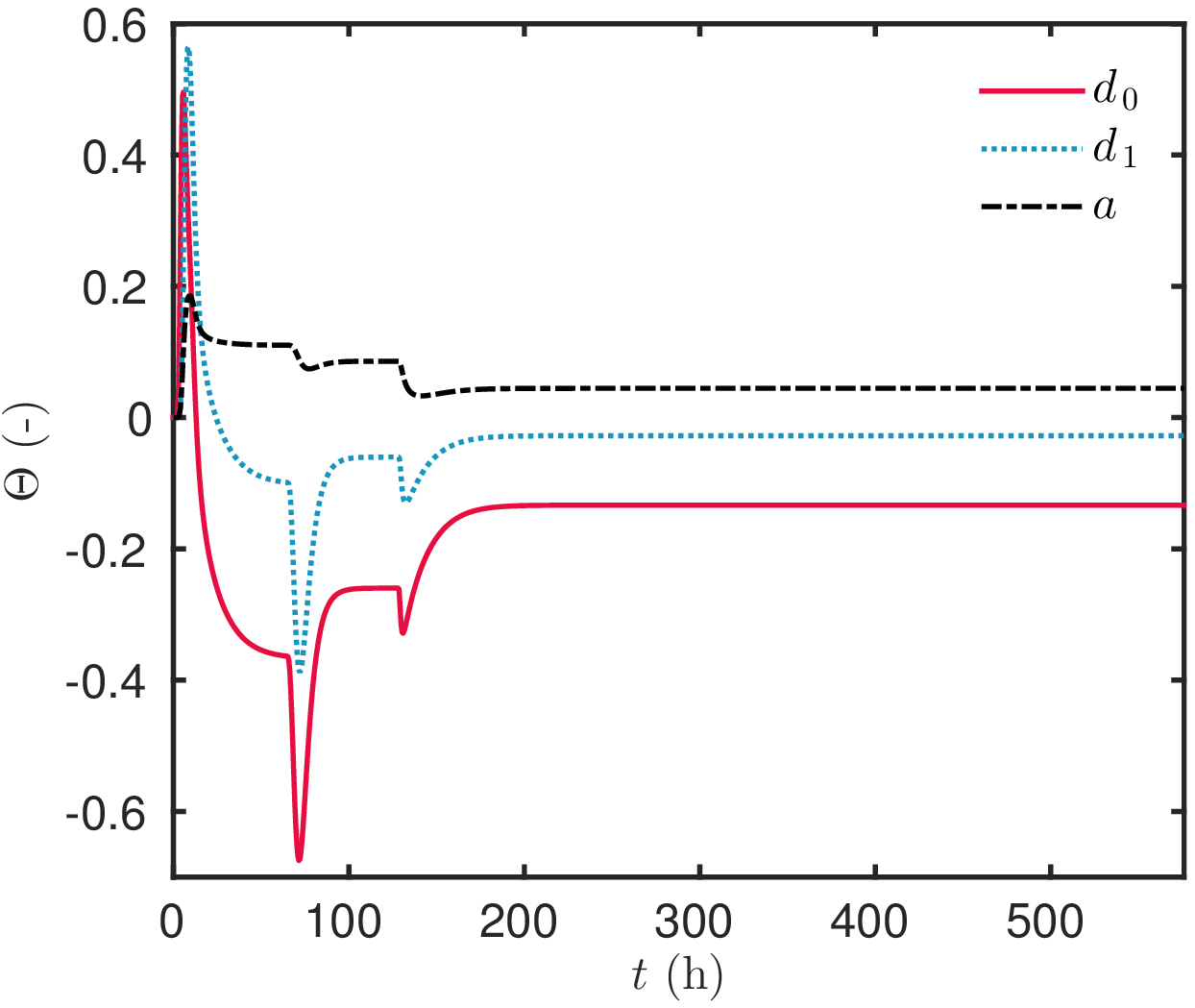}}
\subfigure[][\label{fig:Mstep_S_anti_OED}]{\includegraphics[width=0.48\textwidth]{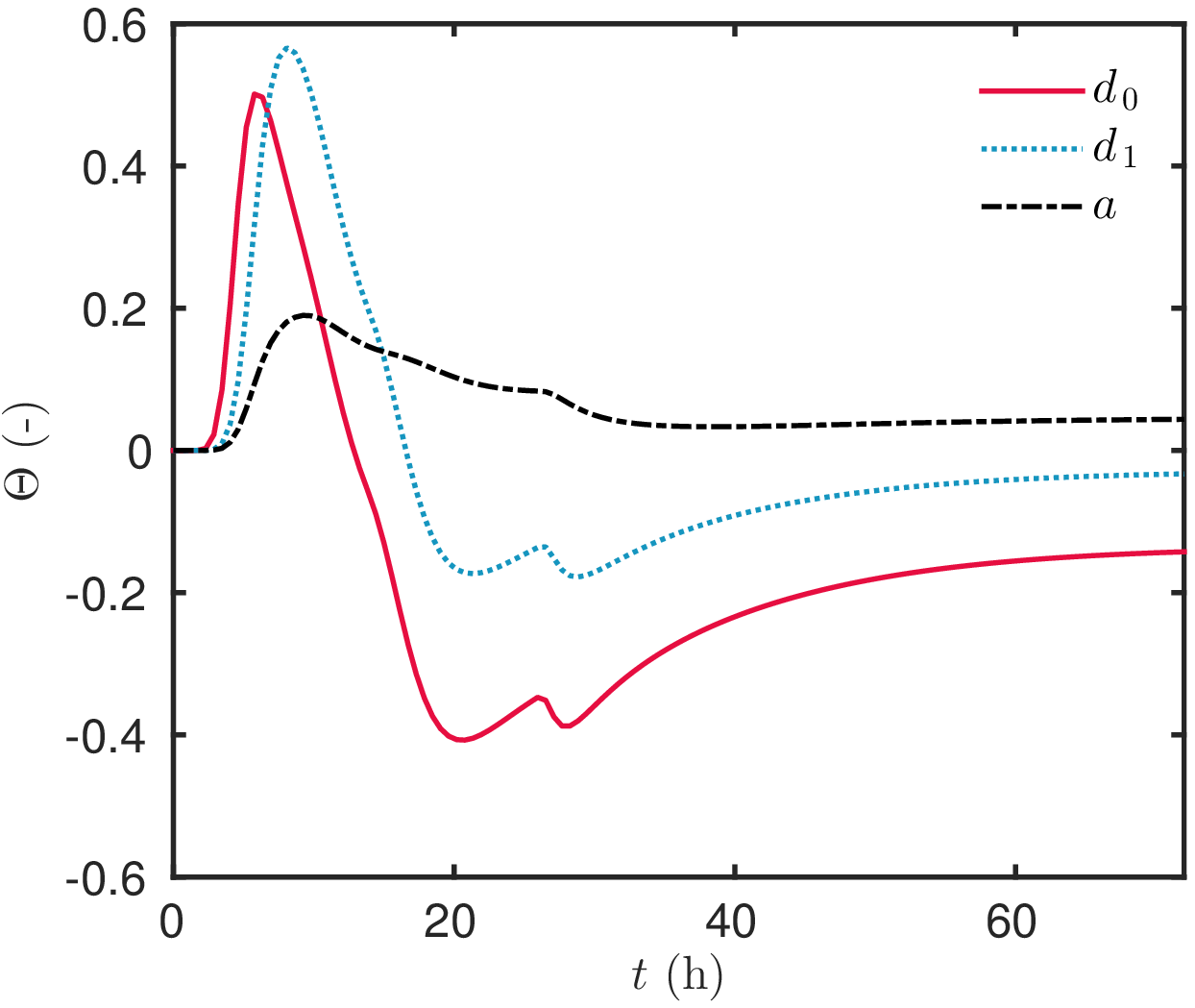}} 
\caption{\small\em Sensitivity coefficients $\Theta$ for parameters $d_{\,0}\,$, $d_{\,1}$ and $a$ for the designs $16$ (OED) (a), $8$ (b) and $1$ (c) ($X \egal X^{\,\circ}$).}
\label{fig:Mstep_S_OED_anti_OED}
\end{center}
\end{figure}

Another type of design is investigated, based on multiple steps of relative humidity. The $16$ possible designs are reported in Table~\ref{tab:Mstep_Design_2parameters}. All the designs have an initial relative humidity $\phi \egal 10 \ \unit{\%}\,$. Two types of consecutive steps are possible: $\phi \egal 33-75-33 \ \unit{\%}$ and $\phi \egal 75-33-75 \ \unit{\%}\,$. The optimal design is sought as a function of the location of the sensor and as a function of the duration of each step of relative humidity. Results are provided in Table~\ref{tab:Mstep_Design_2parameters}. The sensitivity coefficients of $d_{\,0}$ and $d_{\,1}$ are also correlated for these designs as it can be noticed in Figure~\ref{fig:Mstep_S_OED}. Therefore, the estimation of this couple of even three parameters with reasonable accuracy would be a difficult task. Figure~\ref{fig:Mstep_Psi_2param} reveals that the designs with the multiple steps $\phi \egal 75-33-75 \ \unit{\%}$ have a higher magnitude of criterion $\Psi\,$. Moreover, the design $16\,$, corresponding to a duration of each step of $8 \ \unit{days}$, have the highest criterion. The time variation of the sensitivity coefficients are given in Figures~\ref{fig:Mstep_S_OED}, \ref{fig:Mstep_S_8} and \ref{fig:Mstep_S_anti_OED} for the designs $16\,$, $8$ and $1\,$. The coefficients have higher magnitude of variation for design $16$ than for the others, justifying its correspondence to the OED. Figure~\ref{fig:Mstep_Psi_fx_2param} gives the variation of the criterion as a function of the location of the sensors $X^{\,\circ}$ for the design $16\,$. The optimal position is somewhere between $X^{\,\circ} \egal 5 \ \unit{cm}$ and $X^{\,\circ} \egal 6 \ \unit{cm}\,$.

An interesting remark is that the best sensitivity is obtained for an experiment with three steps of $8$ days. For each step, the steady state is reached in the sample. Indeed, it corresponds almost to the duration of the cup experiment for this material, which also carried out measurement until reaching the steady state. This similarity shows a good agreement between the OED results and the cup method.

\begin{table}
\centering
\setlength{\extrarowheight}{.3em}
\begin{tabular}[l]{@{} lcccccccccc}
\hline
\hline
\multirow{2}{*}{Design}  
& Duration of   
& \multicolumn{3}{c}{$(\,d_{\,0} ,\, d_{\,1} \,)$} 
& \multicolumn{3}{c}{$(\,d_{\,0} ,\, a \,)$} 
& \multicolumn{3}{c}{$(\, d_{\,1} ,\, a \,)$} \\
& steps $(\,\unit{days} \,)$
& $\Psi$ & $X^{\,\circ} \ (\,\unit{m}\,)$ & $\mathrm{Corr}$
& $\Psi$ & $X^{\,\circ} \ (\,\unit{m}\,)$ & $\mathrm{Corr}$
& $\Psi$ & $X^{\,\circ} \ (\,\unit{m}\,)$ & $\mathrm{Corr}$ \\
\hline
 \multicolumn{11}{c}{step $ \phi \ (\,\unit{\%} \,)$: 10-33-75-33} \\
1 & 1
& 0.11 & 0.052 & 0.81
& 0.23 & 0.057 & 0.29
& 0.12 & 0.055 & 0.34 \\
2 & 2
& 0.27 & 0.051 & 0.82
& 0.64 & 0.058 & -0.51
& 0.23 & 0.055 & 0.11 \\
3 & 3 
& 0.42 & 0.050 & 0.83
& 1.17 & 0.059 & -0.61
& 0.35 & 0.054 & 0.06 \\
4 & 4 
& 0.54 & 0.050 & 0.83
& 1.63 & 0.059 & -0.66
& 0.44 & 0.053 & -0.10 \\
5 & 5 
& 0.63 & 0.049 & 0.84
& 1.97 & 0.059 & -0.69
& 0.52 & 0.054 & -0.11 \\
6 & 6 
& 0.64 & 0.050 & 0.84
& 1.98 & 0.059 & -0.69
& 0.52 & 0.055 & -0.10 \\
7 & 7 
& 0.65 & 0.051 & 0.85
& 1.99 & 0.059 & -0.69
& 0.53 & 0.054 & -0.10 \\
8 & 8 
& 0.66 & 0.050 & 0.84
& 2.02 & 0.059 & -0.69
& 0.54 & 0.054 & -0.11 \\
\hline
\multicolumn{11}{c}{step $ \phi \ (\,\unit{\%} \,)$: 10-75-33-75} \\
9 & 1
& 1.86 & 0.048 & 0.87
& 3.31 & 0.058 & -0.54
& 1.49 & 0.055 & 0.05 \\
10 & 2
& 5.15 & 0.046 & 0.87
& 8.1  & 0.059 & -0.62
& 3.45 & 0.053 & -0.05 \\
11 & 3 
& 9.12 & 0.047 & 0.88
& 14.4 & 0.060 & -0.65
& 5.73 & 0.053 & -0.10 \\
12 & 4 
& 13.3 & 0.046 & 0.89
& 22.0 & 0.059 & -0.66
& 8.33 & 0.051 & -0.12 \\
13 & 5 
& 16.9 & 0.046 & 0.89
& 27.3 & 0.058 & -0.66
& 10.3 & 0.054 & -0.12 \\
14 & 6 
& 20.0 & 0.045 & 0.90
& 33.5 & 0.059 & -0.66
& 12.0 & 0.052 & -0.13 \\
15 & 7 
& 22.2 & 0.045 & 0.90
& 35.8 & 0.059 & -0.66
& 13.3 & 0.052 & -0.13 \\
\textbf{16} & \textbf{8} 
& \textbf{23.1} & \textbf{0.045} & \textbf{0.90}
& \textbf{36.1} & \textbf{0.059} & \textbf{-0.66}
& \textbf{13.5} & \textbf{0.051} & \textbf{-0.13} \\
\hline
\hline
\end{tabular}
\bigskip
\caption{\small\em Determining the OED (highlighted in bold) for the identification of several parameters, using an experiment with multiple steps of relative humidity.}
\label{tab:Mstep_Design_2parameters}
\end{table}

\subsection{Synthesis on the determination of the OED}

With the given facility, two types of experiments can be carried out with (i) single step and (ii) multiple steps of relative humidity. A total of $20$ designs are possible. The measurement of vapor pressure in the material can be used to estimate the diffusion and advection coefficients $d_{\,0}\,$, $d_{\,1}$ and $a\,$. To provide the best accuracy when solving the parameter estimation problem, the OED has been determined for each design in terms of sensor position and boundary conditions. Two cases has been distinguished. The first one intends to estimate one parameters out of the three. The estimation of the three coefficients at once may lack of accuracy due to the important correlations of the parameters. Thus, the second case aims at estimating two coefficients out of the three. A synthesis of the OED is given in Table~\ref{tab:synthesis_OED}. The choice between the two types of experiments is a compromise between the duration of the total test ($200 \ \unit{h}$ for the single step and $576 \ \unit{h}$ for the multiple steps) and the magnitude of the criterion $\Psi\,$, which is higher for the multiple steps experiments.

\begin{table}
\centering
\setlength{\extrarowheight}{.3em}
\begin{tabular}[l]{@{} c|cc|cc}
\hline
\hline
\multicolumn{5}{c}{Estimation of one parameter} \\
& \multicolumn{2}{c|}{Single step} 
& \multicolumn{2}{c}{Multiple steps} \\
Parameter 
& Boundary condition & Sensor location ($\unit{cm}$) 
& Boundary condition & Sensor location ($\unit{cm}$) \\
$d_{\,0}$ 
& step $ \phi \ (\,\unit{\%} \,)$: 10 -- 75 
& $X \ \in \ \bigl[\, 4 \,,\, 6 \,\bigr]$
& $8 \ \mathsf{days}$ steps, $ \phi \ (\,\unit{\%} \,)$: 10 -- 75 -- 33 -- 75
& $X \ \in \ \bigl[\, 5 \,,\, 6 \,\bigr]$ \\
$d_{\,1}$ 
& step $ \phi \ (\,\unit{\%} \,)$: 10 -- 75 
& $X \ \in \ \bigl[\, 4 \,,\, 6 \,\bigr]$
& $8 \ \mathsf{days}$ steps, $ \phi \ (\,\unit{\%} \,)$: 10 -- 75 -- 33 -- 75
& $X \ \in \ \bigl[\, 3 \,,\, 4 \,\bigr]$ \\
$a$ 
& step $ \phi \ (\,\unit{\%} \,)$: 75 -- 33 
& $X \egal 3$
& $8 \ \mathsf{days}$ steps, $ \phi \ (\,\unit{\%} \,)$: 10 -- 75 -- 33 -- 75
& $X \egal 6$\\
\hline
\multicolumn{5}{c}{Estimation of two parameters} \\
Parameter & Boundary condition & Sensor location 
& Boundary condition & Sensor location ($\unit{cm}$) \\
$(\,d_{\,0} ,\, d_{\,1} \,)$ 
& step $ \phi \ (\,\unit{\%} \,)$: 10 -- 75 
& $X \egal 5$
& $8 \ \mathsf{days}$ steps, $ \phi \ (\,\unit{\%} \,)$: 10 -- 75 -- 33 -- 75
& $X \egal 4.5$ \\
$(\,d_{\,0} ,\, a \,)$
& step $ \phi \ (\,\unit{\%} \,)$: 10 -- 75 
& $X \egal 6$
& $8 \ \mathsf{days}$ steps, $ \phi \ (\,\unit{\%} \,)$: 10 -- 75 -- 33 -- 75
& $X \egal 6$ \\
$(\, d_{\,1} ,\, a \,)$
& step $ \phi \ (\,\unit{\%} \,)$: 10 -- 75 
& $X \egal 5$
& $8 \ \mathsf{days}$ steps, $ \phi \ (\,\unit{\%} \,)$: 10 -- 75 -- 33 -- 75
& $X \egal 5$ \\
\hline
\hline
\end{tabular}
\bigskip
\caption{\small\em Synthesis of the OED for the estimation of one of two parameters, using single or multiple steps experiments.}
\label{tab:synthesis_OED}
\end{table}


\section{Estimation of the parameters}
\label{sec:estim_parameters}

\subsection{Methodology}

As mentioned before, due to the correlation of their sensitivity coefficients, it is uncertain to estimate the couple $\bigl(\, d_{\,0} \,, d_{\,1}\,\bigr)$ with the OED considering single or multiple steps. Therefore, to estimate the three parameters, it is possible to use a couple of parameter with one design (single or multiple steps) and the third one with the other. It can be noticed that the criterion $\Psi$ is higher for the OED with multiple steps of relative humidity than a single step. Moreover, for this OED, the sensitivity coefficients of the parameters $\bigl(\, d_{\,1},\, a\,\bigr)$ have lower correlation that for the parameters $\bigl(\, d_{\,0} \,, a\,\bigr)\,$. Therefore, to identify three parameters $\bigl(\, d_{\,0} \,, d_{\,1} \,, a\,\bigr)\,$, the design $16$ with multiple steps $10-75-33-75 \ \unit{\%}$ will be used to estimate the couple $\bigl(\, d_{\,1},\, a\,\bigr)$ and the design $2$ with the single step $10-75 \ \unit{\%}$ for the remaining parameter $d_{\,0}\,$. The sensor was located according to the OED.

We denote as  $\bigl(\,  d_{\,0}^{\,\circ} \,, d_{\,1}^{\,\circ} \,, a^{\,\circ} \,\bigr)$ and $\bigl(\,  d_{\,0}^{\,\mathrm{apr}} \,, d_{\,1}^{\,\mathrm{apr}} \,, a^{\,\mathrm{apr}} \,\bigr)$ the estimated and \emph{a priori} parameters, respectively. Table~\ref{tab:unkown_parameters} gives the numerical values of the later. To estimate the parameter, the following cost function is defined, in its matrix form, using the Euclidean norm, as: 
\begin{align*}
\mathrm{J} \, \bigl(\,d_{\,0} \,, d_{\,1} \,, a \,\bigr) 
& \egal \biggl(\, u_{\,\mathrm{exp}}^{\,1} \moins u_{\,1} \bigl(\, X^{\,\circ} \,, t \,, d_{\,0} \,\bigr) \, \biggr)^T \cdot
\biggl(\, u_{\,\mathrm{exp}}^{\,1} \moins u_{\,1} \bigl(\, X^{\,\circ} \,, t \,, d_{\,0} \,\bigr) \, \biggr) \nonumber \\
& \plus
\biggl(\, u_{\,\mathrm{exp}}^{\,2} \moins u_{\,2} \bigl(\, X^{\,\circ} \,, t \,, d_{\,1} \,, a \,\bigr) \, \biggr)^T \cdot
\biggl(\, u_{\,\mathrm{exp}}^{\,2} \moins u_{\,2} \bigl(\, X^{\,\circ} \,, t \,, d_{\,1} \,, a \,\bigr) \, \biggr) \,,
\end{align*}
where subscripts $1$ and $2$ denote the design with single and multiple steps of relative humidity, respectively. The quantity $u_{\,\mathrm{exp}}$ stands for the measurement resulting from the experimental facility and the OED. It is evaluated on a discrete grid of time $t_{\,\mathrm{exp}}\,$. The quantity $u$ is the solution of the direct problem --- Eqs.~\eqref{eq:HAM_equation2}--\eqref{eq:HAM_ic} ---, considering the boundary and initial conditions of the design.

The cost function $\mathrm{J}$ is minimised using function \texttt{fmincon} in the \texttt{Matlab} environment, providing an efficient interior-point algorithm with constraint on the unknown parameters \cite{Byrd2000}. Here, the constraints are defined with upper and lower bound for the parameters:
\begin{align*}
d_{\,0}^{\,\circ} \ \in \ \bigl[ \, 0.01 \,, 10 \, \bigl] \cdot d_{\,0}  \,, &&
d_{\,1}^{\,\circ} \ \in \ \bigl[ \, 0.01 \,, 10 \, \bigl] \cdot d_{\,1} \,, &&
a^{\,\circ} \ \in \ \bigl[ \, 0.001 \,, 10 \, \bigl] \cdot a \,.
\end{align*}
As mentioned in Section~\ref{sec:exp_facility}, the \emph{a priori} advection parameter has been estimated without experiments. Thus a larger interval was provided as a constraint for this parameter.


\subsection{Results and discussion}

\begin{figure}
\begin{center}
\subfigure[][\label{fig:phi_1step}]{\includegraphics[width=0.48\textwidth]{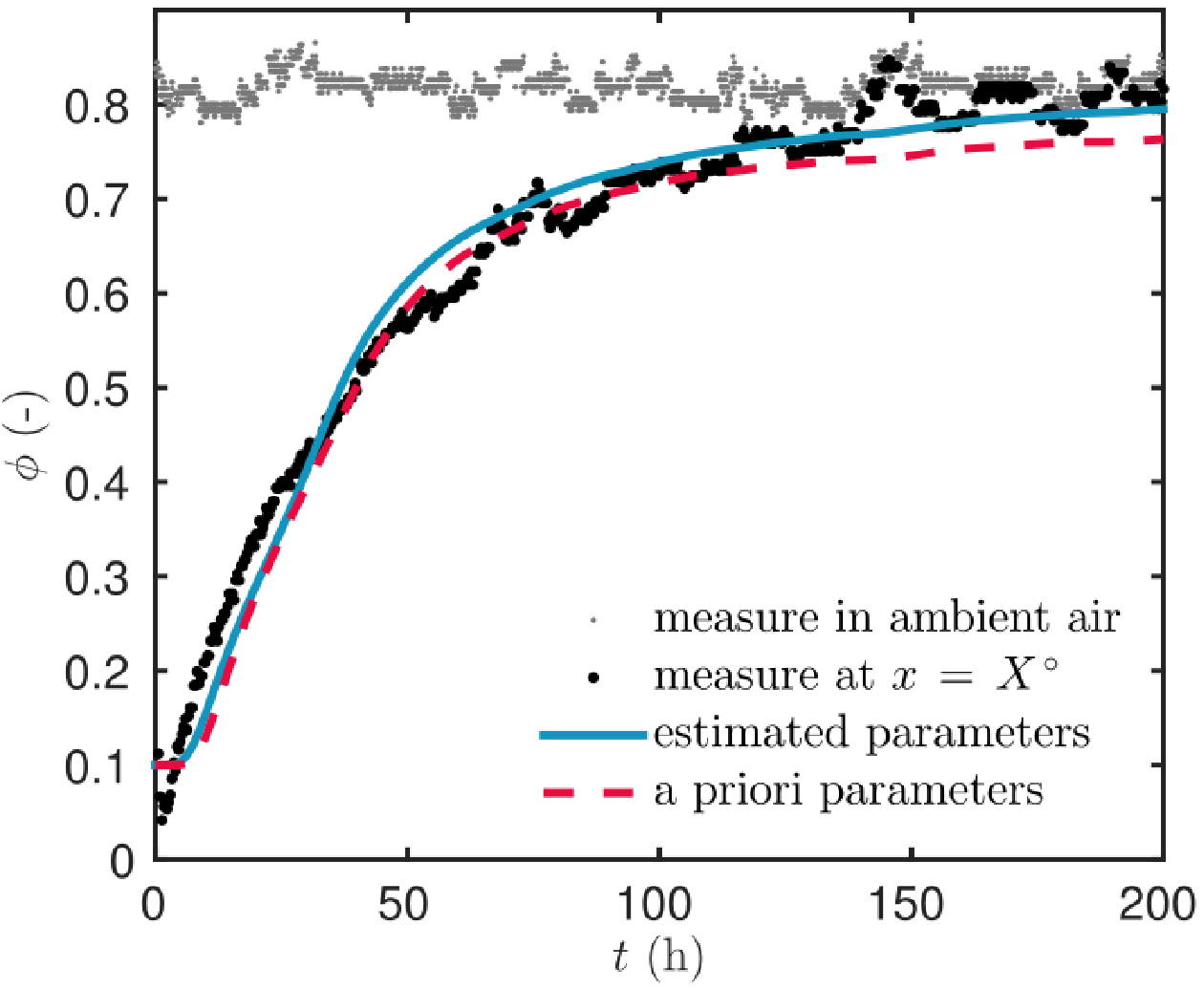}}
\subfigure[][\label{fig:phi_Mstep}]{\includegraphics[width=0.48\textwidth]{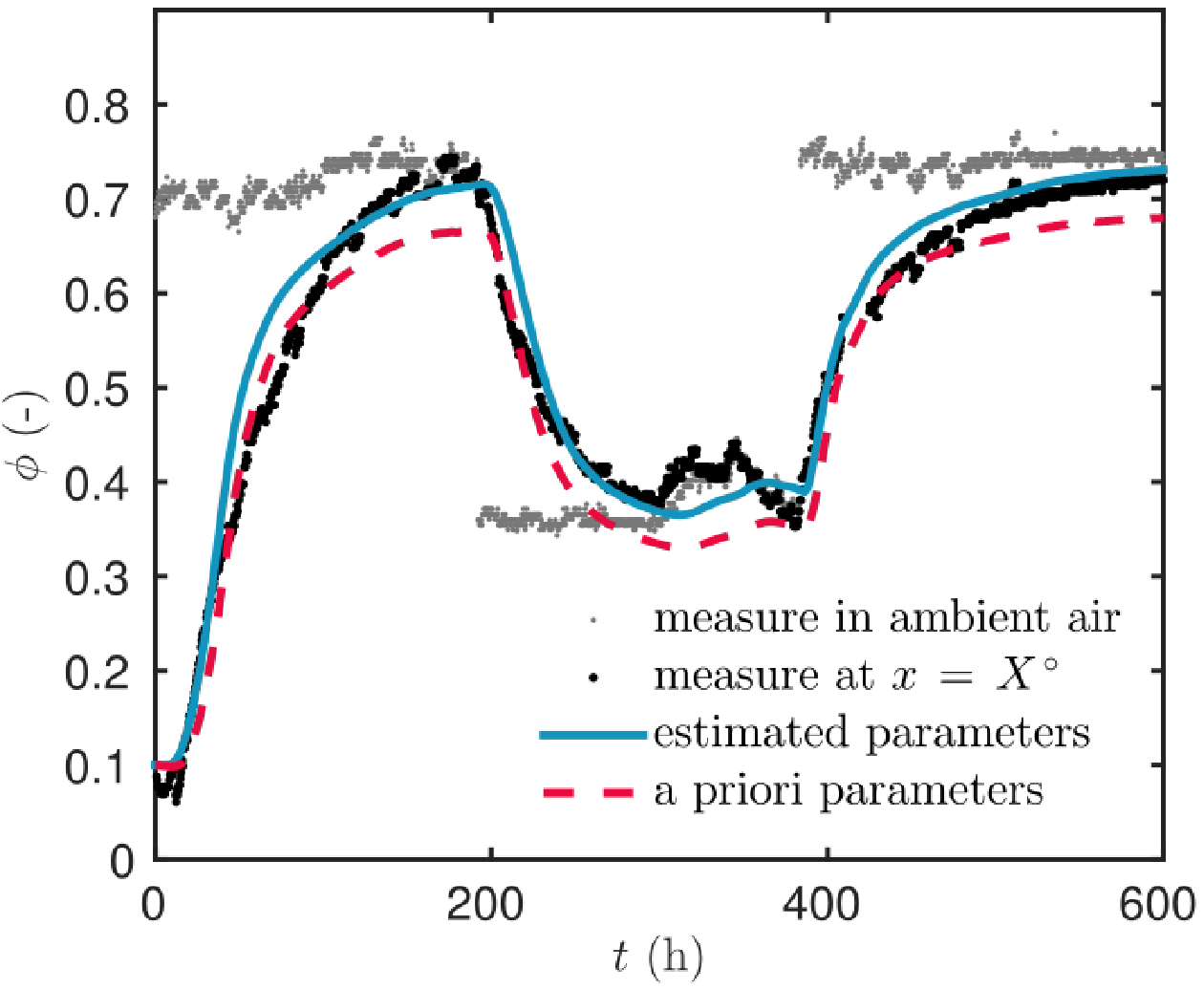}}
\caption{\small\em Comparison of the experimental data and the numerical results.}
\label{fig:est_phi}
\end{center}
\end{figure}

Results of the estimated parameters are reported in Table~\ref{tab:unkown_parameters}. It took only $54$ computations of the direct model for both experiments to minimise the cost function $\mathrm{J}\,$. In addition, a comparison between the experimental measured data and the results of the numerical models with \emph{a priori} and estimated parameters is given in Figure~\ref{fig:est_phi}. It can be noted that the regulation by salt solutions in the climatic chamber is not very precise. For the design with a single step, the ambient relative humidity is around $\phi_{\,1} \egal 80 \ \unit{\%}$ instead of $75 \ \unit{\%}\,$.

The three estimated parameters have higher values than the \emph{a priori} ones. However, the values are not so different to consider that the methodology used for the search of the OED is valid. The estimated parameters for the vapour permeability is almost two times higher. The water diffusion resistance factor $\mu$ of the material, corresponding to the ratio between the vapour permeability and the air permeability, varies between $\mu \egal 4.5$ and  $\mu \egal 2.8$ for dry and wet states. As this factor remains lower than $1$, one can say that the estimated vapour permeability is physically acceptable. It can be mentioned that in previous study \cite{Rouchier2017}, with similar experiments for this material and a physical model including only diffusion transfer, authors estimated a resistance factor higher than $1\,$. They concluded that an \emph{apparent} permeability, even not-physical, was estimated and that other physical phenomena should be taken into account. The present study goes further by including the moisture advection, providing a more complete and reasonable physical model. As noticed in Table~\ref{tab:unkown_parameters}, a ratio of $10$ is noticed between the \emph{a priori} and estimated advection parameters $a\,$, highlighting that the advection in the material is less important. The mass average velocity in the material is of the order of $0.01 \ \unit{mm/s}\,$. Moreover, the moisture permeability and advection coefficients have the same order of magnitude, confirming the results from \cite{Abahri2016}, highlighting the importance of both advection and diffusion transfer in porous materials.

These results are also correlated with the ones from \cite{Duforestel2015} where authors questioned the measurement of the vapour permeability according to the standard \textsf{ISO 12572}. The method appears not to take into account the variation of the total pressure within the cup. Thus, advective vapour transfer is overlooked. In the case of highly air permeable materials, it implies that the vapour permeability may be underestimated.

\begin{figure}
\begin{center}
\subfigure[][\label{fig:Res_1step}]{\includegraphics[width=0.48\textwidth]{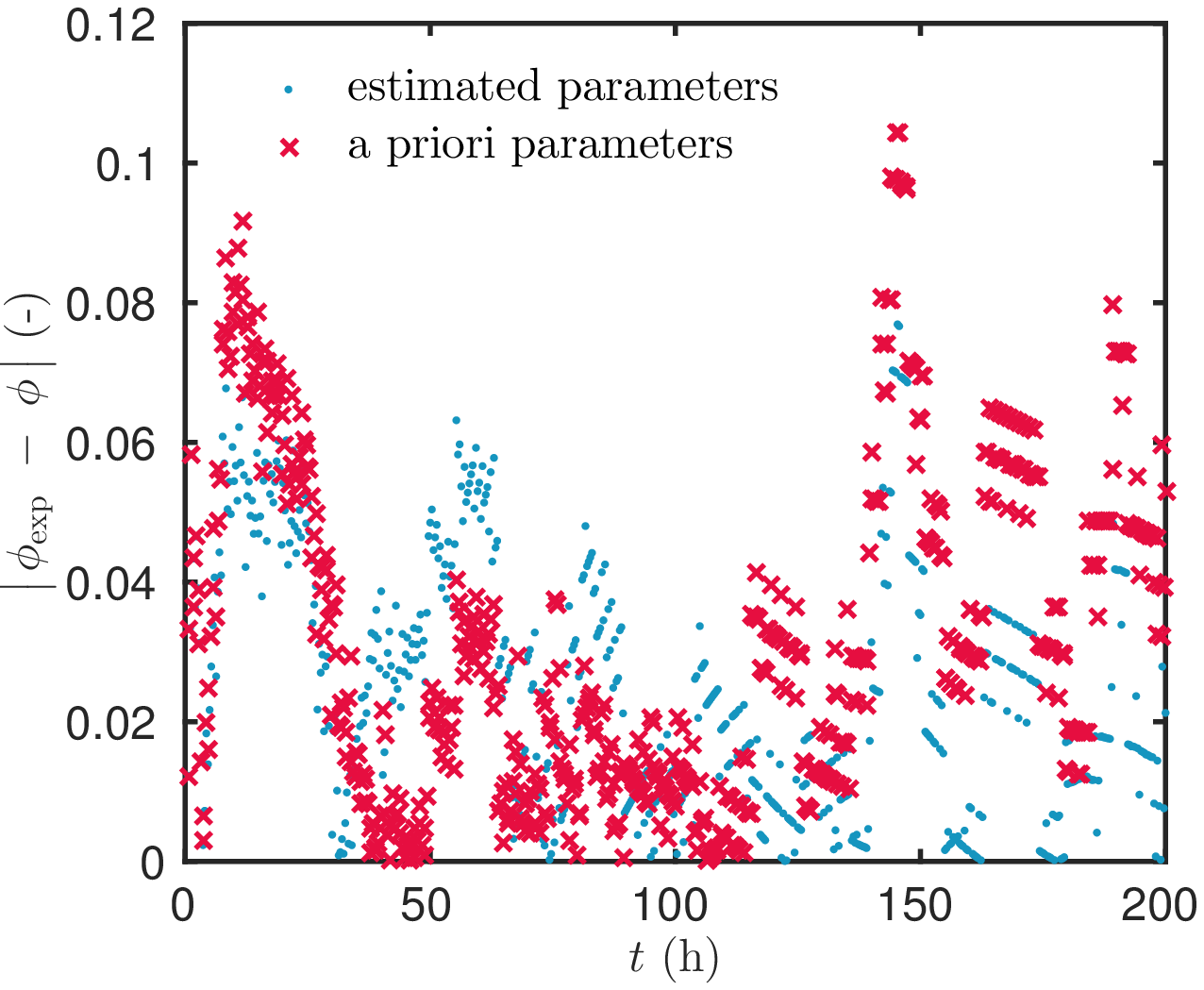}}
\subfigure[][\label{fig:Res_Mstep}]{\includegraphics[width=0.48\textwidth]{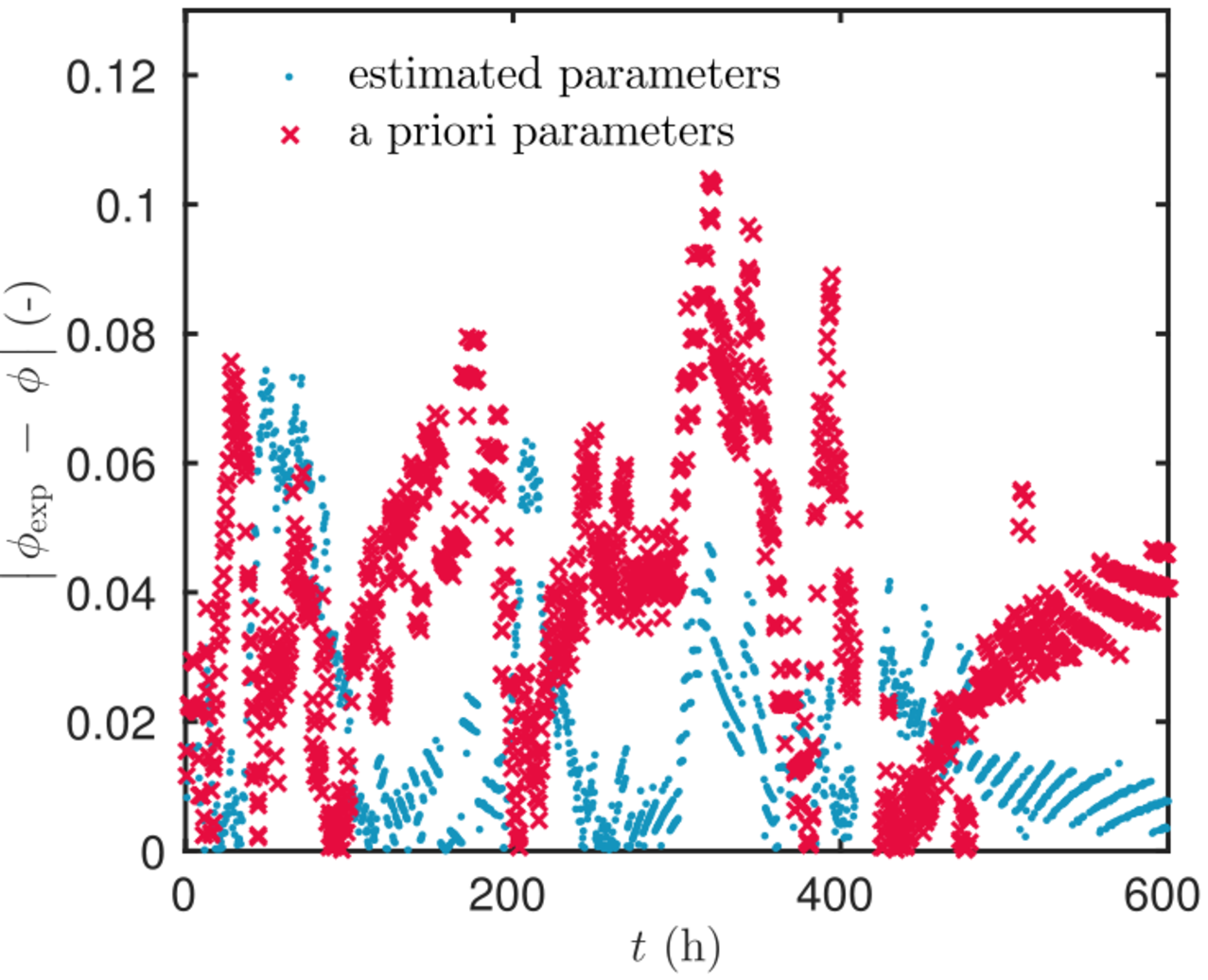}}
\caption{\small\em Residual between the experimental data and the numerical results.}
\label{fig:residu}
\end{center}
\end{figure}

Figure~\ref{fig:residu} gives the residual between the numerical results and the experimental data for both experiments. The residual is lower for the results with the estimated parameters. It is also uncorrelated, highlighting a satisfying estimation of the parameters. Nevertheless, it can be noted in Figure~\ref{fig:est_phi} that some discrepancies remain between the experimental data and the results of the numerical model. This can be specifically observed at $t \egal 150 \ \unit{h}\,$, in Figure~\ref{fig:phi_1step}, and at $t \egal 350 \ \unit{h}\,$, in Figure~\ref{fig:phi_Mstep}. These discrepancies might be explained by two reasons. First, it was assumed that the vapour permeability varies as a first-degree polynomial of the relative humidity  in Eq.~\eqref{eq:mat_hypothesis} and that the advection parameter is independent of the fields. To fit better with the experiments, this assumption may be revisited by increasing the degree of the vapour permeability polynomial and considering a field dependent non-constant advection parameter. Another alternative could be to search for time varying coefficients as it was done in \cite{Berger2016b}. The second reason is that only moisture transfer is included in the physical model, neglecting the variation of the temperature in the samples. Even if the experimental design is carried out under isothermal conditions, the temperature varies slightly in the material due to phase change and coupling effects between heat and moisture transfer. Moreover, temperature should influence at least the advection coefficient $a$. Thus, the physical model should be improved for achieving a better fitting with the experimental data.

\begin{figure}
\begin{center}
\subfigure[][\label{fig:pdfd0_ft}]{\includegraphics[width=0.48\textwidth]{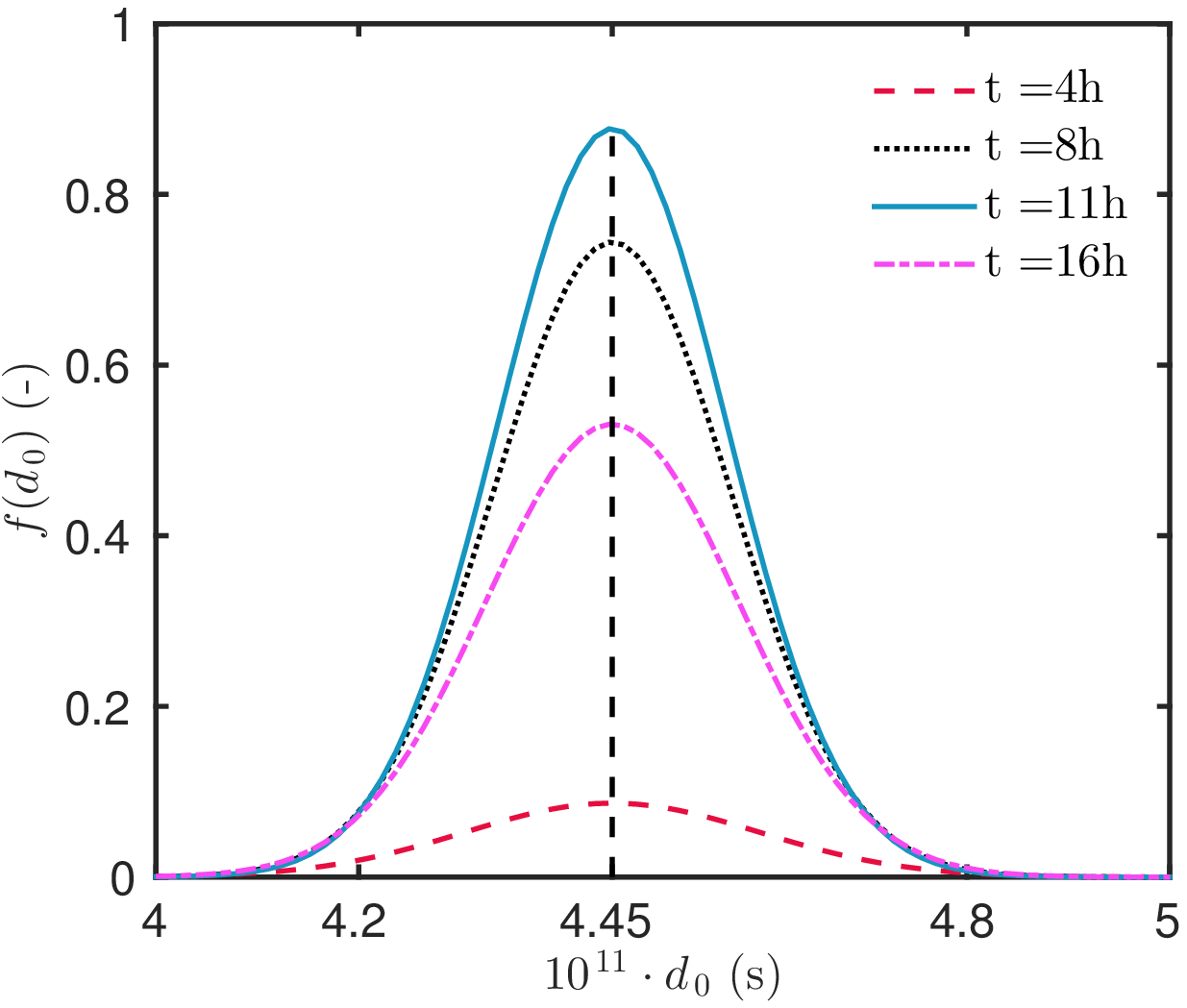}}
\subfigure[][\label{fig:pdfd1_ft}]{\includegraphics[width=0.48\textwidth]{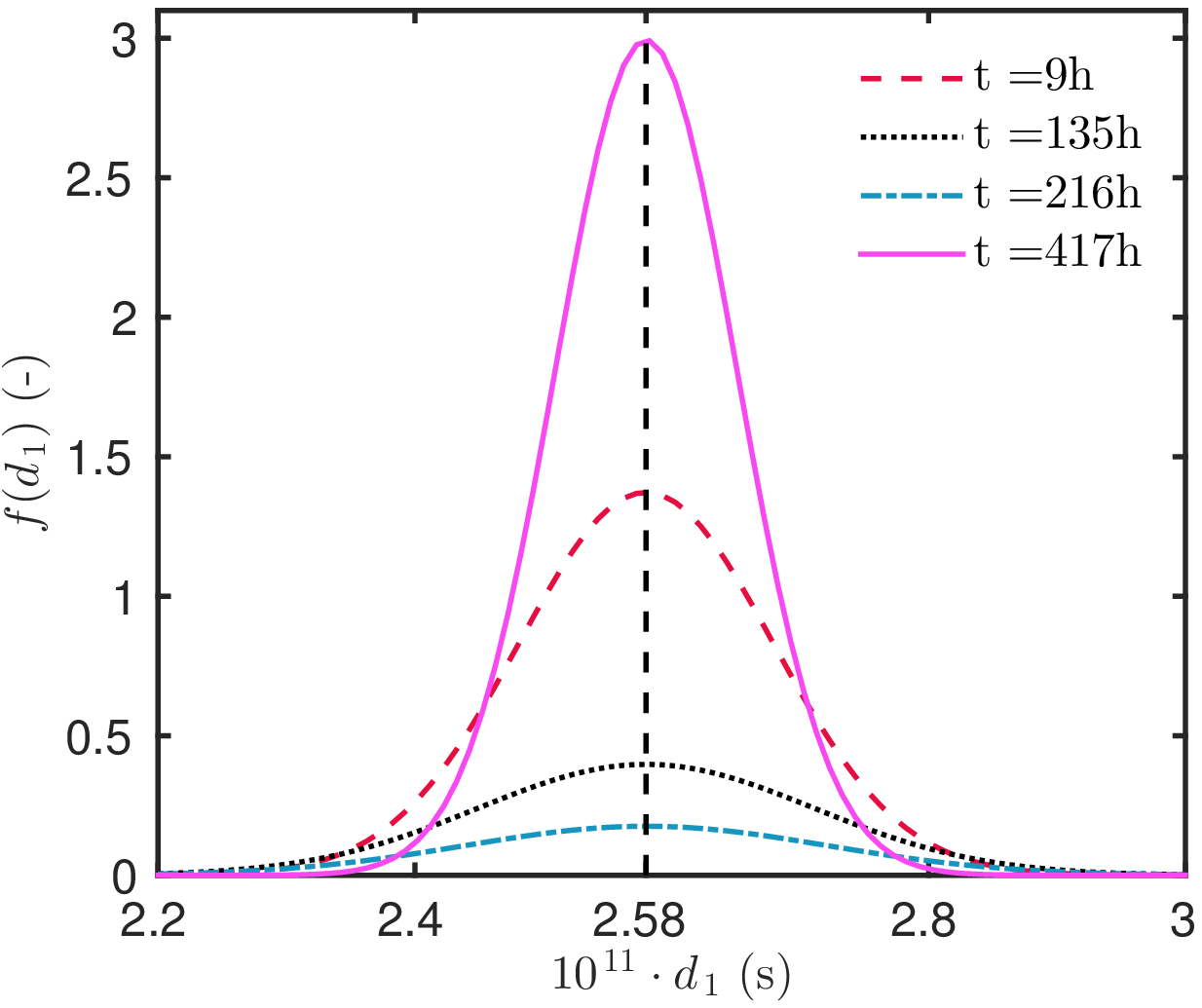}}
\subfigure[][\label{fig:pdfa_ft}]{\includegraphics[width=0.48\textwidth]{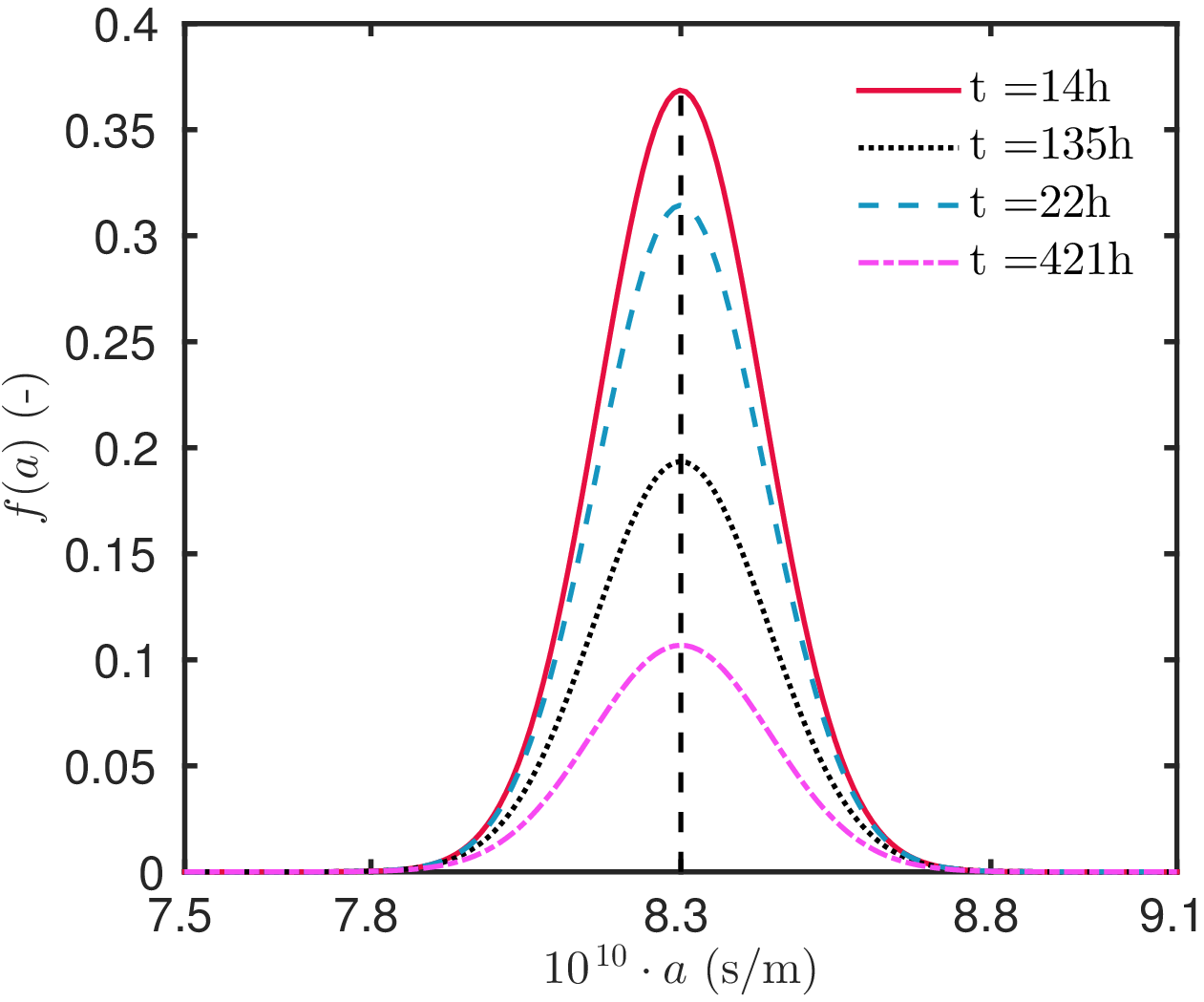}}
\caption{\small\em Probability density function approximated for the estimated parameters.}
\label{fig:pdf}
\end{center}
\end{figure}

As mentioned in Section~\ref{subsec:OED}, using Eqs~\eqref{eq:probability_pm_1} and \eqref{eq:probability_pm_2}, the local probability density functions of the estimated parameters can be computed using the sensitivity functions of the parameters. A normal distribution with a standard deviation of measurement $\sigma_{\,\phi} \egal 0.15 \unit{\%}$ is assumed. As the field is varying in time, the probability of the parameter and its standard deviation are also depending on time as illustrated in Figure~\ref{fig:pdf}. In Figure~\ref{fig:1step_S_OED}, the sensitivity function of this parameter reaches its maximum for $t \egal 11\; \unit{h}\,$. It also corresponds to the time where the probability function is maximum, as noted in Figure~\ref{fig:pdfd0_ft}. A precise estimation of the parameter is obtained at the moment when the sensitivity function reaches its maximum. Similar observations can be done when comparing the sensitivity function of parameters $d_{\,1}$ and $a$ in Figures~\ref{fig:Mstep_S_OED} and the probability in Figures~\ref{fig:pdfd1_ft} and \ref{fig:pdfa_ft}. The interesting point is to get an information on the local probability of the estimated parameter, at a low computational cost. It is computed directly from the sensitivity function. It is not necessary to carry out numerical-based statistical inferences that have a non-negligible computational cost as mentioned in \cite{Berger2016b} for instance.

\begin{table}
\centering
\setlength{\extrarowheight}{.3em}
\begin{tabular}[l]{@{} lcc }
\hline
\hline
Parameters & \emph{a priori} value & estimated value \\
Moisture permeability $\bigl[\, \unit{s} \,\bigr]$ 
& $d_{\,0}^{\,\mathrm{apr}} \egal 2.33 \e{-11} $ 
& $d_{\,0}^{\,\circ} \egal 4.45 \e{-11} $ \\
Moisture permeability $\bigl[\, \unit{s} \,\bigr]$ 
& $d_{\,1}^{\,\mathrm{apr}} \egal 5.68 \e{-11} $ 
& $d_{\,1}^{\,\circ} \egal 2.58 \e{-11} $ \\
Advection parameter $\bigl[\, \unit{s/m} \, \bigr] $
& $a^{\,\mathrm{apr}} \egal 7.2 \e{-11} $ 
& $a^{\,\circ} \egal 8.3 \e{-10} $ \\
\hline
\hline
\end{tabular}
\bigskip
\caption{\small\em \emph{A priori} and estimated parameters.}
\label{tab:unkown_parameters}
\end{table}


\section{Conclusion}

This paper presented a study conducted for the estimation of moisture-dependent material properties using \emph{in-situ} measurements, with two mains objectives. The first one was to compare the results to previous studies using a similar facility and material \cite{Rouchier2017}. The second goal was to propose a more detailed physical model improved by adding a moisture advection term, commonly neglected in the literature. The unknown parameters to be determined were the moisture diffusion and advection coefficients of a wood fibre material. The primer one was assumed to vary as a first-degree polynomial of the relative humidity, while the latter was defined as a constant parameter. Thus, the inverse problem involved three parameters  $\bigl(\, d_{\,0}, \, d_{\,1},\, a \,\bigr)$ of a diffusive-convective moisture transfer problem.

The estimation of the unknown parameters has been accomplished using an existing set-up based on two chambers with prescribed relative humidity and temperature. The set-up enables to submit samples to single or multiple steps of relative humidity under isothermal conditions. Before performing the experiments, the OED was searched in terms of optimal boundary conditions and location of sensors within the material. These conditions ensure to provide the highest accuracy of the identification method. This approach was described in \cite{Berger2017} and applied to this study for the existing facility among $20$ possible designs. It has been carried out considering \emph{a priori} values of the unknown parameters determined with standard methods \cite{Rafidiarison2015}. For the estimation of one parameter among the three, results have shown that a single step of relative humidity from $\phi \egal 10 \ \unit{\%}$ to $\phi \egal 75 \ \unit{\%} \,$, with a sensor located between $4$ and $6 \ \unit{cm}$ is the optimal design for the estimation of the moisture permeability coefficients. For the advection coefficient, a single step of $\phi \egal 75 \ \unit{\%}$ to $\phi \egal 33 \ \unit{\%} \,$, with a sensor located at $ X^{\,\circ} \egal 3 \ \unit{cm}\,$, is the optimal design. It has also been demonstrated that the sensitivity coefficients of parameters $d_{\,0}$ and $d_{\,1}$ have strong correlations for those experiments. Therefore, it would be a difficult task to identify all the parameters with only one experiment. Nevertheless, the OED was also determined for the identification of coupled parameters. A sample submitted to  multiple relative humidity steps of $\phi \egal 10-75-33-75 \ \unit{\%}\,$, with a sensor placed at $ X^{\,\circ} \egal 5 \ \unit{cm} \,$ appears as the best option. It is important to mention that this investigation is valid for the \emph{a priori} values of material under investigation as well for the conditions of the facility (length of the samples, boundary conditions, \etc).

After these investigations, two OEDs have been used to perform two different experiments. The first one, a single step of relative humidity from $\phi \egal 10 \ \unit{\%}$ to $\phi \egal 75 \ \unit{\%}\,$, is used to estimate the parameter $d_{\,0}$. The second experiment is based on multiple steps of relative humidity $\phi \egal 10-75-33-75 \ \unit{\%}$ to determine the couple $\bigl(\, d_{\,1}\,, a \,\bigr)\,$. For each experiment, the sensors have been placed according to the OED. Using these measures, a cost function has been defined and minimised using an interior-point state-of-the-art algorithm. The estimated moisture permeability coefficients are twice higher than the \emph{a priori} values of the parameters, being a physically acceptable moisture diffusion resistance factor. The advection parameter has been determined, corresponding to the air speed inside the material of the order of $\mathsf{v} \egal 0.01 \ \unit{mm/s}\,$, which is in accordance with the results from other studies \cite{Berger2017, Belleudy2016}.

To conclude, the values of the estimated parameters are close to the ones of \emph{a priori} parameters used for the search of the OED, validating the approach. Moreover, thanks to the computation of the sensitivity function within the OED methodology, a probability distribution function of the estimated parameters can be approximated at a lower computational cost. As highlighted in Section~\ref{sec:estim_parameters}, some discrepancies still remain between the experiments and the numerical calibrated model, which motivates further research to be focused on upgrading the physical model by considering the combined heat and moisture transfer to improve the calibration. In addition, the assumptions on the vapour permeability and advection coefficients may be revisited. An interesting improvement of the methodology would be to consider both fields of relative humidity and average moisture content measurements when solving the parameter estimation problem. In that case, the sensitivity coefficient for both fields should be computed  and the new results of the OED methodology should be then analyzed.


\bigskip

\subsection*{Acknowledgments}
\addcontentsline{toc}{subsection}{Acknowledgments}

This work was partly funded by the French Environment and Energy Management Agency (ADEME), the ``Assembl\'ee des pays de Savoie'' (APS) and the French National Research Agency (ANR) through its Sustainable Cities and Buildings program (MOBAIR project ANR-12-VBDU-0009).

\bigskip


\appendix
\section*{Nomenclature}

\begin{tabular*}{0.7\textwidth}{@{\extracolsep{\fill}} | c  l l| }
\hline
\multicolumn{3}{|c|}{\emph{Latin letters}} \\
$a$ & advection coefficient & $[\mathsf{s/m}]$ \\
$d$ & moisture diffusion & $[\mathsf{s}]$ \\
$c$ & moisture storage capacity & $[\mathsf{kg/m^3/Pa}]$ \\
$h$ & vapour convective transfer coefficient & $[\mathsf{s/m}]$ \\
$k$ & permeability & $[\mathsf{s}]$ \\
$L$ & length & $[\mathsf{m}]$ \\
$\Pc$ & capillary pressure & $[\mathsf{Pa}]$ \\
$\Ps$ & saturation pressure & $[\mathsf{Pa}]$ \\
$\Pv$ & vapour pressure & $[\mathsf{Pa}]$ \\
$R_{\,v}$ & water gas constant & $[\mathsf{J/kg/K}]$\\
$T$ & temperature & $[\mathsf{K}]$ \\
$\mathsf{v}$ & air velocity & $[\mathsf{m/s}]$ \\
\multicolumn{3}{|c|}{\emph{Greek letters}} \\
$\phi$ & relative humidity & $[-]$ \\
$\rho$ & specific mass & $[\mathsf{kg/m^3}]$ \\
\hline
\end{tabular*}


\bigskip
\addcontentsline{toc}{section}{References}
\bibliographystyle{abbrv}

\bigskip

\end{document}